\documentclass[11pt]{article}
\usepackage{graphicx}
\usepackage{amsmath,latexsym,amsbsy,amssymb, color}
\usepackage{graphicx}
\font\bigbf=cmbx10 at 16pt

\def\ds{\displaystyle}
\def\forall{\hbox{for all}~}

\def\ve{\varepsilon}

\def\J{{\cal J}}
\def\I{{\cal I}}

\def\R{{\mathbb R} }%I\!\!R}
\def\bfC{{\bf C}}

\def\implies{\Longrightarrow}

\def\prob{\hbox{Prob.}}

\def\vs{\vskip 2em}
\def\vsk{\vskip 4em}
\def\v{\vskip 1em}

\def\cR{{\cal R}}

\def\S{{\cal S}}

\def\cN{{\cal N}}

\def\Tilde{\widetilde}
\def\Hat{\widehat}
\def\bega{\begin{array}}
\def\enda{\end{array}}
\def\begi{\begin{itemize}}
\def\endi{\end{itemize}}
\def\meas{\hbox{meas}}

\def\P{\hbox{Prob.}}
\def\bel{\begin{equation}\label}
\def\eeq{\end{equation}}
\def\sqr#1#2{\vbox{\hrule height .#2pt
\hbox{\vrule width .#2pt height #1pt \kern #1pt
\vrule width .#2pt}\hrule height .#2pt }}
\def\square{\sqr74}
\def\endproof{\hphantom{MM}\hfill\llap{$\square$}\goodbreak}

\newtheorem{theorem}{Theorem}[section]

\newtheorem{lemma}{Lemma}[section]
\newtheorem{proposition}{Proposition}[section]
\newtheorem{remark}{Remark}[section]
\newtheorem{definition}{Definition}[section]
\newtheorem{example}{Example}[section]

\begin{document}

\title{\bigbf Markovian Solutions to Discontinuous ODEs}
\vs
\author{Alberto Bressan$^{(1)}$, Marco Mazzola$^{(2)}$, and Khai T. Nguyen$^{(3)}$\\ 
\\
 {\small $^{(1)}$ Department of Mathematics, Penn State University, }\\  
 {\small $^{(2)}$ IMJ-PRG, CNRS, Sorbonne Universit\'e,}\\
  {\small $^{(3)}$ Department of Mathematics, North Carolina State University. }\\ 
 \\  {\small e-mails: ~axb62@psu.edu, ~marco.mazzola@imj-prg.fr,  ~ khai@math.ncsu.edu}
 }
\maketitle

\begin{abstract}
Given a possibly discontinuous, bounded  function $f:\R\mapsto\R$,
we consider the set of generalized flows, obtained by assigning a probability measure on the set of Carath\'eodory solutions to the ODE ~$\dot x = f(x)$.
The paper provides a complete characterization of all such flows which 
have a Markov property in time. This is achieved in terms of (i) a 
positive, atomless
measure supported on the set $f^{-1}(0)$ where $f$ vanishes, 
(ii) a countable number of Poisson random variables, determining the  waiting times
 at points in $f^{-1}(0)$, and (iii) a countable set of numbers $\theta_k\in [0,1]$,
 describing the probability of moving up or down, at isolated points where
two distinct trajectories can originate. 
\end{abstract}

\vsk

\section{Introduction}
\label{s:1}
\setcounter{equation}{0}

Consider a scalar  ODE with possibly
discontinuous  right hand side:
\bel{ode}
\dot x~=~f(x).\eeq
Given an initial condition
\bel{id}x(0)~=~x_0\,,\eeq
we recall that $t\mapsto x(t)$ is a {\em Carath\'eodory solution} 
of the Cauchy problem 
(\ref{ode})-(\ref{id}) if
\bel{CS}
x(t)~=~x_0+\int_0^tf(x(s))\, ds\qquad\qquad \forall t\geq 0\,.
\eeq
When the function $f$ is not Lipschitz continuous, it is well known that 
this Cauchy problem can admit multiple solutions.   
Because of this non-uniqueness, 
in \cite{E} it was proposed to study ``generalized flows", described by a probability
measure on the set of all Carath\'eodory solutions.   See also \cite{FL} for a
related approach.
 %Given an initial point $x_0$ and $t\geq 0$, one thus consider a probability 
%distribution on the set of solutions
In this direction,  the main goal of the present paper is to describe all stochastic flows 
compatible with the ODE, which have the Markov property. Toward this goal, in the first part of the paper we describe all 
possible deterministic semigroups:

\begin{definition}\label{d:DS} By a {\bf deterministic semigroup} compatible with the
ODE (\ref{ode}) we mean  a map $S:\R\times \R_+\mapsto \R$,
with the properties
\begi\item[(i)] $S_t(S_s(x_0))~=~S_{t+s}(x_0)$, ~$S_0(x_0)~=~ x_0$.
\item[(ii)] For each $x_0\in \R$, the map $t\mapsto S_t x_0$
is a solution to the Cauchy problem (\ref{ode})-(\ref{id}).
\endi
\end{definition}
Notice that here we do not require any continuity w.r.t.~the initial point $x_0$. Throughout the following, as basic assumptions we consider:
\begi
\item[{\bf (A1)}] {\it The function $f:\R\mapsto\R$ is  bounded and regulated, i.e.~it admits left and right limits $f(x-)$, $f(x+)$ at every point $x$.   }
\item[{\bf (A2)}] {\it If $y$ is  a point  where
either $f(y-)\cdot f(y+)=0$ or else $f(y-)> 0> f(y+)$, then $f(y)=0$.}
\endi
\begin{remark} {\rm The assumption   {\bf (A2)} corresponds to the ``no-jamming"
condition used in \cite{Binding}. 
It is needed to 
rule out examples such as
$$\dot x~=~f(x)~=~\left\{\bega{rl} 1\quad &\hbox{if}\quad x<0,\cr
-1\quad &\hbox{if}\quad x>0,\cr
a\quad &\hbox{if}\quad x=0,\enda\right.\qquad\qquad x(0)=0,$$
which has no Carath\'eodory solution if $a\not= 0$.}
\end{remark}
\begin{remark} {\rm 
It is well known that every regulated function has at most countably many points of discontinuity.  Namely, the set
\bel{Df}
D_f~=~\Big\{x\in\R\,;~f(x-)\neq f(x+)~~\mathrm{or}~~f(x-)=f(x+)\neq f(x)\Big\}
\eeq
is at most countable. We observe that, under the assumptions {\bf (A1)-(A2)}, the set
of zeroes
\bel{zset}f^{-1}(0)~=~\bigl\{x\in\R\,;~f(x)=0\bigr\}\eeq
is closed.
Indeed, assume that $f(x_n)=0$ for a strictly increasing sequence $x_n\to \bar x$.  
By {\bf (A1)} this implies that the left limit exists and satisfies
$f(\bar x-) =0$. Hence {\bf (A2)}  yields $f(\bar x)=0$. 
Similarly, if $x_n$ decreases to $\bar x$, then $f(\bar x+)=0$ and hence
{\bf (A2)}  again implies $f(\bar x)=0$.
}
\end{remark}
In Section~\ref{s:3}, we will  show that, in order to uniquely determine
a deterministic semigroup compatible with (\ref{ode}),
three additional ingredients are needed:
\begi
\item[{\bf (Q1)}]
A continuum (i.e., atomless) positive measure $\mu$ supported on 
the set of zeroes (\ref{zset}) of $f$.
Strictly increasing  trajectories $t\mapsto x(t)= S_t(x_0)$ of the semigroup are then defined by the identity
$$t~=~\int_{x_0}^{x(t)} {dy\over f(y)}+ \mu\bigl([x_0, x(t)]\bigr),$$
while a similar formula holds for decreasing ones.
Notice that, if $f^{-1}(0)$ is countable, then necessarily $\mu=0$.
\item[{\bf (Q2)}] A countable set of points $\S\subseteq f^{-1}(0)$, where the dynamics
is forced to stop.   Among the (possibly many) solutions of (\ref{ode})
starting from  $x_0\in \S$, we thus choose $S_t (x_0)=x_0$.
\item[{\bf (Q3)}] A map $\Phi:\Omega^*\mapsto \{-1,1\}$, defined
on a set $\Omega^*\subset\R$ of isolated points from where
both an increasing and a decreasing solution of (\ref{ode}) can originate.
For $x_0\in \Omega^*$, setting $\Phi(x_0)= 1$ selects the increasing
solution, while $\Phi(x_0)=-1$ selects the decreasing one.
\endi
In the second part of the paper we study Markov processes, whose sample
paths are solutions to (\ref{ode}).
Given a deterministic semigroup compatible with (\ref{ode}), 
in Section~\ref{s:5} we show that
a Markov semigroup can be obtained by adding two more ingredients to the 
list {\bf (Q1)--(Q3)}. Namely:
\begi
\item[{\bf (Q4)}]  A countable set $\S^*\subset f^{-1}(0)$ and a
map $\Lambda: \S^*\mapsto \,]0, +\infty[\,$, describing the random
waiting time of a trajectory which reaches a point $y_j\in\S^*$.
More precisely, we assume that a solution initially at $y_j\in \S^*$ remains at 
$y_j$ for a random waiting time $Y_j\geq 0$, then starts moving.
All these random times will be mutually independent, with Poisson distribution:
\bel{Tk}\prob\{Y_j > s\} ~=~e^{-\lambda_j s}\qquad\mathrm{with}\qquad \lambda_j = \Lambda(y_j).\eeq
%with $\lambda_j = \Lambda(y_j)$.   
%For notational convenience, we 
%also allow the cases $\lambda_k=0$,
%so that $T_k= +\infty$, and $\lambda_k=+\infty$, so that $T_k= 0$ with probability one.
\v
\item[{\bf (Q5)}] A map $\Theta:\Omega^*\mapsto [0,1]$, defined 
on the countable set $\Omega^*\subset\R$ of points from which
both an increasing and a decreasing solution of (\ref{ode}) can originate.
For $x_k\in \Omega^*$, the value $\theta_k =\Theta(x_k)$ gives the probability that,
when the solution starting from $x_k$ begins to move, it will be increasing.
Of course, $1-\theta_k$ is then the probability that the solution will be decreasing.
\endi
Conversely,  in Section~\ref{s:6} we prove that every Markov semigroup
whose sample paths are solutions to (\ref{ode}) can be obtained 
by first constructing a deterministic semigroup $S$  as in {\bf (Q1)--(Q3)},
and then adding the random waiting times $Y_j$ in {\bf (Q4)} and the upward vs.~downward
probabilities $\theta_k$ in {\bf (Q5)}.

\begin{example} {\rm Given $0<\alpha<1$, consider the 
ODE
\bel{xa}\dot x~=~{|x|^\alpha\over 1-\alpha}.\eeq
This provides a classical example of an ODE with continuous right hand side, with
multiple solutions starting at  the origin.   This 
ODE is compatible with  
two deterministic semigroups $S,\Tilde S:\R\times\R_+\mapsto\R$. 
The first one satisfies 
$S_t(0) = 0$ for all $t\geq 0$, while the second one satisfies
$\Tilde S_t(0)= \ds t^{1\over 1-\alpha}$.
In addition, there is a one-parameter family of Markovian flows.
Indeed, for each  $0<\lambda<\infty$, one can consider the 
random time $T$ with Poisson distribution
$$\prob\{ T>s\}~=~\ds e^{-\lambda s}.$$
A random solution starting at $x=0$ 
 remains at the origin up to 
time $T$, then it is strictly increasing for $t>T$, namely $x(t) = (t-T)^{1\over 1-\alpha}$.
Notice that, by taking the limits $\lambda\to 0$ or $\lambda\to +\infty$, we recover the 
deterministic semigroups $S$ and $\Tilde S$, respectively.
}
\end{example}

The remainder of the paper is organized as follows.   In Section~\ref{s:2}
we briefly review some basic properties of 
scalar ODEs with discontinuous right hand side.
For the main results on existence and uniqueness of solutions
we refer to  \cite{Binding}. These are based on the elementary solution formula
\bel{sol} \int_{x_0}^{x(t)} {1\over f(y)}\, dy~=~t\,,\eeq
valid on intervals where $f$ has a.e.~the same sign, and $1/f$ is integrable.
 
Section~\ref{s:3} is concerned with deterministic semigroups.
Theorem~\ref{car_sem}
describes the most general semigroup $S:\R\times \R_+\mapsto \R$ compatible with the ODE (\ref{ode}). 
In Section~\ref{s:5} we construct a family of Markov semigroups
whose sample paths are solutions to (\ref{ode}).  Each of these semigroups
is determined by the additional data introduced at {\bf (Q1)--(Q5)}.
Finally, Theorem~\ref{t:61}  in Section~\ref{s:6} shows that every Markov semigroup compatible with (\ref{ode}) is one of the above family.

In the companion paper \cite{BMN}, the authors show that every deterministic
semigroup compatible with (\ref{ode}) can be obtained as the pointwise limit of the flows 
generated by a sequence of ODEs $\dot x=f_n(x)$ with smooth right hand sides.
Moreover, every Markovian semigroup can be obtained as
limit of a sequence of diffusion processes with smooth drifts and with diffusion coefficients approaching zero.

Ordinary differential equations with discontinuous right hand side
have been studied for their own interest \cite{Binding, B88, BC,
Fil, Sentis}, as well as for
their several applications to feedback control \cite{AB, CLSW, MP, PS}.
They also play a key role in the theory  of transport equations \cite{ABC, AC,CR,DPL}, and in the 
 analysis of solutions to hyperbolic conservation laws  by the method of characteristics \cite{BS, CM, D1, D2}.    For a 
survey of numerical methods we refer to \cite{DL}.

\section{Review of scalar discontinuous ODEs}
\setcounter{equation}{0}
\label{s:2}
As a preliminary, we collect here
 some basic results on the existence and properties of solutions to 
a  Cauchy problem with possibly discontinuous right hand side.
In the following theorem, 
the first two statements can already be found in \cite{Binding}, where a more general setting was considered.
However,  the closure of the solution set strongly relies on the assumption
{\bf (A1)} that the function $f$ is regulated.
\v
\begin{theorem}\label{t:1} Consider the Cauchy problem (\ref{ode})-(\ref{id}),
assuming that {\bf (A1)-(A2)} hold. Then 
\begi
\item[(i)]  For every $x_0\in\R$ there exists at least one Carath\'eodory  solution, defined for all times $t\geq 0$. 
\item[(ii)]
Every solution is monotone (either increasing or decreasing).
\item[(iii)] The set of all solutions is nonempty and closed
w.r.t.~the topology of uniform convergence on bounded sets.
\endi
\end{theorem}

{\bf Proof.} 
{\bf 1.} We begin by proving  that (\ref{ode})-(\ref{id}) admits a 
Carath\'eodory solution. 
By the assumption {\bf (A1)} there exists a constant $M>0$ such that
\bel{M}
|f(x)|~\leq~M\qquad\forall x\in \R.
\eeq
Let any $x_0\in \R$ be given.  If $f(x_0) = 0$ then trivially
 $x(t)\equiv x_0$ is a  solution of (\ref{ode}). 
On the other hand, if  $f(x_0) \not=  0$, then by  the assumption {\bf (A2)} we must have
\[
\hbox{either}~~f(x_0+)~>~0\quad\mathrm{or}\quad f(x_0-)~<~0\,.
\]

To fix the ideas, assume $f(x_0+)>0$,
 the other case being entirely similar. 
We can then find $\delta,h>0$ such that 
$$h~\leq~f(x)~\leq~M\qquad\qquad\hbox{for all }~x\in \,]x_0, x_0+\delta].$$
The map
$$x~\mapsto~ t(x)~\doteq~\int_{x_0}^x {dy\over f(y)}$$
is thus strictly  increasing in the interval $[x_0,x_0+\delta]$. 
Indeed, for the above inequality, for all $ x_0\leq x_1<x_2\leq x_0+\delta$ it holds
\[
{x_2-x_1\over M}~\leq~t(x_2)-t(x_1)~=~\int_{x_1}^{x_2}{dy\over f(y)}~\leq~{x_2-x_1\over h}\,.
\]
%$$t(x_2)-t(x_1)~=~\int_{x_1}^{x_2}{dy\over f(y)}~\leq~{x_2-x_1\over h}\quad\forall x_0\leq x_1<x_2\leq x_0+\delta\,.$$
This implies that the inverse map 
$t\mapsto x(t)$ is strictly increasing and Lipschitz continuous.
Hence it is  absolutely continuous, and satisfies
\[
\dot x(t)~=~\left({d\over dx} t(x)\right)^{-1}~=~f(x(t))\quad \hbox{for a.e.}~~ t\in [t(x_0),t(x_0+\delta)].
\]
Thus, $x(t)$ is a solution of (\ref{ode}), defined on the interval $[t(x_0),t(x_0+\delta)]$.

We claim that this solution can be extended for all $t\in [0, +\infty[\,$.
Indeed, if $[0,T[\,$ is a maximal domain where the solution can be constructed,
since $f$ is bounded by $M$, the limit $x(T)\doteq \lim_{t\to T-} x(t)$ exists.
By the previous arguments, the solution can then be extended
to $[0, T+\delta]$ for some $\delta>0$, contradicting the maximality assumption.
\medskip

{\bf 2.} To prove that every  solution is monotone, let $t\mapsto x(t)$ be a 
Carath\'eodory  solution of (\ref{ode}). Assume, on the contrary, that $x(\cdot)$ is 
neither monotone increasing nor decreasing. Without loss of generality 
we can assume that 
\[
x(t_1)~=~x(t_2)~<~x(\bar{t})\qquad\hbox{for some}~~t_1<\bar{t}<t_2\,.
\]
Since $x(\cdot)$ is absolutely continuous, there exists $\tau_1\in \,]t_1,t_2[\,$ such that  $x(\cdot)$ is differentiable at $\tau_1$ and
$$\dot x(\tau_1)>0,\qquad x(t_1)~<~x(\tau_1)~<~x(\bar{t})\,.$$  
Therefore there exists $\delta_0>0$ such that
\[
x(\tau_1)\,<\, x(\tau_1+s)\qquad\ \forall ~s\in \,]0,\delta_0[\,.
\]
Since $f$ has  a right limit  at $x(\tau_1)$, we deduce
\bel{ineq1}
%f(x(\tau_1)-)~=~\lim_{h\to 0+}{1\over h} 
%\int_{-h}^0 f(x(\tau_1+s))~ds~=~
0~<~\dot x(\tau_1)~=~\lim_{h\to 0+}{1\over h} \int_{0}^{h}f(x(\tau_1+s))~ds~=~f(x(\tau_1)+).
\eeq
On the other hand, since $x(t_2)<x(\tau_1)<x(\bar{t})$, we consider the time
$$\tau_2~\doteq~ \min\bigl\{ t\geq \bar t~~:~~x(t) = x(\tau_1)\bigr\}.$$
The above definition implies 
$x(t)> x(\tau_1)$ for all $t\in [\bar t, ~\tau_2[\,$ while 
$$x(\tau_2)~=~\lim_{t\to \tau_2-} x(t)~=~x(\tau_1).$$
As a consequence, we have
\[
f(x(\tau_1)+)~=~\lim_{h\to 0+}{1\over h}
\int_{-h}^0f(x(\tau_2+s))~ds~=~\lim_{h\to 0+}{x(\tau_2)-x(\tau_2-h)\over h}~\leq~0,
\]
reaching a contradiction with (\ref{ineq1}).
\medskip

{\bf 3.} To prove (iii), consider  a sequence  of solutions $x_n(\cdot) $  of (\ref{ode}), converging  to a function $\bar{x}(\cdot) $ uniformly on bounded sets. We claim that $\bar x$ is also a Carath\'eodory  solution. 
To fix the ideas,
assume that all solutions  $x_n$ are monotone increasing. Hence 
$\bar{x}$ is monotone increasing as well.  We need to show that 
\bel{ke1}
\bar{x}(b)~=~\bar{x}(a)+\int_{a}^bf(\bar{x}(s))~ds
\eeq
for any interval $[a,b]$.  Toward this goal, consider the set of times
\[\bega{rl}
I_0&\ds=~\Big\{t\in [a,b]~;\quad \hbox{there exists $\ve>0$ such that}\\[3mm] 
&\qquad\qquad \hbox{
$\bar x(\cdot)$ is constant either on $[t,t+\ve]$ or on $[t-\ve,t]$}\Big\}.
\enda
\]
We  claim that 
\bel{cla1}
\bar x(t)~\in~f^{-1}(\{0\})\qquad\forall t\in I_0\,.
\eeq
Indeed, by  the definition of $I_0$, for any $t_0\in I_0$ there exist $t_1<t_2$
with $t_0\in [t_1,t_2]$   such that 
\[
\bar x(t)~=~\bar  x(t_0)\qquad\forall t\in [t_1,t_2]\,.
\]
Hence,
\[
\lim_{n\to \infty}~\int_{t_1}^{t_2}f(x_n(t))~dt~=~\lim_{n\to \infty} \bigl[ x_n(t_2) - x_n(t_1)\bigr] ~=~0\,.
\]
Since $x_n$ is monotone increasing, we have $f(x_n(s))= x_n'(s) \geq 0$
for a.e.~$s$.  This implies
\[
\lim_{n\to\infty}~\|f(x_n)\|_{{\bf L}^1([t_1,t_2])}~=~0\,.
\]
By possibly taking a subsequence, this 
yields the pointwise convergence
$f(x_n(t))\to 0$ for a.e.~$t\in [t_1,t_2]$.

As a consequence, one of the three following cases must hold:
\begi
\item[(i)]  $f(x(t_0))~=~0$,
\item[(ii)]  $f(x(t_0)-)~=~0$,
\item[(iii)]  $f(x(t_0)+)~=~0$.
\endi
By the assumption {\bf (A2)}, both (ii) and (iii) imply  (i), proving our claim (\ref{cla1}).
In particular, this implies
\bel{lim1}\lim_{n\to\infty} f(x_n(t))~=~f(\bar x(t))~=~0
\qquad\hbox{ for a.e. ~} t\in I_0\,.\eeq
Next, from the definition of $I_0$ it follows that  the function $\bar{x}(\cdot)$ is one-to-one (strictly increasing) on $[a,b]\backslash I_0$. In particular, this implies that the set 
\[
\bigl\{t\in [a,b]\setminus I_0\,;~~\bar{x}(t)\in D_f\bigr\}
\]
is at most countable.
Since $f$ is continuous outside $D_f$ and $x_n(t)$ converges to $\bar{x}(t)$ for every $t\in [a,b]$, we have
\bel{lim2}
\lim_{n\to\infty}~f(x_n(t))~=~f(\bar x(t))\qquad \hbox {for a.e.~}~
t\in [a,b]\setminus I_0\,.
\eeq
By (\ref{lim1})-(\ref{lim2}), the dominated convergence theorem yields
\[
\lim_{n\to\infty}~\int_{a}^bf(x_n(t))~ds~=~\int_a^bf(\bar x(t))~dt\,.
\]
This proves (\ref{ke1}).
\endproof

\begin{example}\label{e:21} {\rm To appreciate the importance of  the assumption {\bf (A1)}, 
consider the following Cauchy problem, proposed in \cite{Binding}.
\bel{CP3} \dot x~=~f(x)~=~\left\{\bega{rl} 0\quad&\hbox{if}\quad x= 1, {1\over 2}, {1\over 3},\ldots
\\[3mm]
1\quad &\hbox{otherwise,}\enda\right.\qquad\qquad x(0)=0.\eeq
Then for every $n\geq 1$ the function
$$x_n(t)~=~\min\Big\{ {1\over n}, \,t\Big\}$$
is a Carath\'eodory solution, but the uniform limit  $x(t) = \ds\lim_{n\to \infty} x_n(t)~=~0$
does not satisfy (\ref{CP3}).   Notice that the function $f$ in (\ref{CP3}) is not regulated.
}
\end{example}

\section{Semigroups compatible with the ODE}
\label{s:3}
\setcounter{equation}{0}

In general, the Cauchy problem
for the discontinuous ODE (\ref{ode})-(\ref{id}) admits several different solutions.  
Our present goal is to study all possible semigroups $S:\R\times\R_+\mapsto\R$
compatible with the ODE (\ref{ode}), 
in the sense of Definition~\ref{d:DS}.
%Our eventual  goal is to put probability measures on the 
%set of solutions.  More precisely, we wish to describe the family of 
%all Markov processes whose sample paths are solutions to (\ref{ode}),
%with probability one.

We first introduce a few definitions.
For any function $f$, denote
$$f^+(x)\doteq\max\{f(x),0\},\qquad\qquad f^-(x)= \min\{f(x),0\}.$$
If $f$ is regulated, so are $f^+, f^-$.
Moreover, we define the sets
%{\color{blue} 
%Consider the  set of points accessible in finite time from the left
%and from the right, respectively:
%\bel{A-}
%\A^-~\doteq~
%\left\{x\in \R\,;~~~\int_{x-\ve}^x {1\over f^+(y)}\, dy
%\,<\,+\infty\quad\hbox{for some}~\ve>0\right\},\eeq
%\bel{A+}
%\A^+~\doteq~
%\left\{x\in \R\,;~~~\int_x^{x+\ve} {1\over f^-(y)}\, dy
%\,>\,-\infty\quad\hbox{for some}~\ve>0\right\}.\eeq
%In other words, $\bar x\in \A^-$ if and only if there exists a strictly increasing
%Carath\'eodory solution
%to (\ref{ode}), say $x:[-\tau, 0]\mapsto\R$, with $\tau>0$,  such that $x(0)=\bar x$.
%Similarly, $\bar x\in \A^+$ if there exists a strictly decreasing Carath\'eodory solution to 
%(\ref{ode}) which reaches $\bar x$ in finite time.
%In addition we define the set of points reaching out to the left and to the right,
%respectively:}
\bel{R-}
\cR^-~\doteq~
\left\{x\in \R\,;~~~\int_{x-\ve}^x {1\over f^-(y)}\, dy
\,>\,-\infty\quad\hbox{for some}~\ve>0\right\},\eeq
\bel{R+}
\cR^+~\doteq~
\left\{x\in \R\,;~~~\int_x^{x+\ve} {1\over f^+(y)}\, dy
\,<\,+\infty\quad\hbox{for some}~\ve>0\right\}.\eeq
Notice that $\bar x\in \cR^-$ if there exists a strictly decreasing
Carath\'eodory solution
to (\ref{ode}), say $x:[0,\tau]\mapsto\R$, with $\tau>0$, such that $x(0)=\bar x$.
Similarly, $\bar x\in \cR^+$ if there exists a strictly increasing Carath\'eodory solution to 
(\ref{ode}) which starts at $\bar x$.

From the above definitions, it immediately follows that $\cR^-$ is
open to the left and $\cR^+$ is open to the right.   Namely
$$x\in \cR^-\qquad\implies\qquad ]x-\delta, x]\subset\cR^-
\qquad\hbox{for some} ~~\delta>0,$$
$$x\in \cR^+\qquad\implies\qquad [x, x+ \delta[\subset\cR^+
\qquad\hbox{for some} ~~\delta>0.$$
For future use, we also define the sets
\bel{O*}
\Omega^*~\doteq~\cR^-\cap\cR^+,\qquad \Omega_0~\doteq~\R\setminus\bigl( \cR^-\cup\cR^+).
\eeq
%\bel{O0}
%\Omega_0~\doteq~\R\setminus\bigl( \cR^-\cup\cR^+).\eeq
Notice that $\Omega^*$ is the set of points from which two distinct 
solutions of (\ref{ode}) can initiate:
one strictly decreasing, and one strictly increasing.
On the other hand, $\Omega_0$ is a set of points from which neither an increasing
nor a decreasing solution can initiate.    If $x_0\in \Omega_0$, the unique solution
to the Cauchy problem (\ref{ode})-(\ref{id}) is the constant one:
$x(t)\,=\, x_0$ for all $t\geq 0$.   
\begin{remark}{\rm 
By  {\bf (A2)}, the set of zeroes $f^{-1}(0)$ is closed. Moreover,
the assumptions {\bf (A1)-(A2)}  imply the inclusion 
$\Omega_0\subseteq f^{-1}(0)$.
We also observe that  the set  $\Omega^*$ consists of  isolated points. 
Therefore, there can be at most countably many  such points.
}\end{remark}

\begin{remark}\label{r:32}{\rm In general, the set 
\bel{smz}f^{-1}(0) \cap (\cR^+\cup\cR^-)\eeq
has measure zero, hence its interior is empty.   However, it may well be uncountable.
For example let ${\bf C}\subset [0,1]$ be the Cantor set, which we write in the form
\bel{Cantor}
{\bf C}~=~[0,1]\,\backslash\,\bigcup_{n=0}^{\infty}~
\bigcup_{k=1}^{2^{n}}~I^n_{k}
\eeq
where $I_{k}^{n}=\left]a^n_k, b^n_k\right[$ are disjoint open intervals with length $3^{-(n+1)}$. We define 
\[
f(x)~=~\left\{ \bega{cl}
0\qquad &\hbox{if}~~ x\in \bfC\,\cup \,]-\infty,0]\,\cup \,[1,+\infty[\,,\\[3mm]
\min\left\{\big|x-a^n_k\big|^{1/3},\big|x-b_k^n\big|^{1/3}\right\}\qquad &\hbox{if}~~
x\,\in\, I_k^n%~=~] a_k^n\,,\,b_k^n[ 
 \,. \enda\right.
\]
It is easy to check that $f$ is continuous and 
\begin{eqnarray*}
\int_{[0,1]\setminus \bfC}{1\over f(x)}~dx&=&\sum_{n=0}^{\infty}\sum_{k=1}^{2^n}\int_{I_k^n}{1\over f(x)}dx~=~3\sum_{n=0}^{\infty}\sum_{k=1}^{2^n}\left({3^{-(n+1)}\over 2}\right)^{2/ 3}\cr\cr
&=&{3^{1/3}\over 2^{2/3}}\cdot\sum_{n=0}^{\infty}\left(2\over 3^{2/3}\right)^n~<~+\infty.
\end{eqnarray*}
This provides an example where 
the set  $f^{-1}(0)\cap \cR^+= {\bf C}\backslash \{1\}$ is uncountable.
As a consequence, this set can support a nontrivial atomless measure $\mu$.}
\end{remark}

%\subsection{Constructing a semigroup compatible with the ODE.}
We now consider the ODE (\ref{ode}), together
with a positive, atomless measure $\mu$ supported on $f^{-1}(0)$, 
a countable set $\S\subseteq f^{-1}(0)$, and a map $\Phi:\Omega^*\mapsto 
\{-1, 1\}$, as described at  {\bf (Q1)--(Q3)} in the Introduction.  
We claim that these assignments uniquely determine a semigroup $S:\R\times\R_+\mapsto\R$, compatible with the ODE (\ref{ode}).

\begin{definition}\label{d:increase} We say that an open interval $]a,b[$ is a 
{\bf domain of increase}
if
\bel{din}
]a,b[~\cap~\S~=~\emptyset,\qquad \mu([c,d])+\int_c^d {dx\over f^+(x)}~<~+\infty
\qquad\forall [c,d]\subset\,]a,b[\,.\eeq
Similarly, we say that  $]a,b[$ is a 
{\bf domain of decrease}
if
\bel{dde}
]a,b[~\cap~\S~=~\emptyset,\qquad \mu([c,d])-\int_c^d {dx\over f^-(x)}~<~+\infty
\qquad\forall [c,d]\subset\,]a,b[\,.\eeq
\end{definition}

If $]a,b[$ and $]a', b'[$ are two intervals of increase having non-empty
intersection, then their union $]a,b[\,\cup\,]a', b'[$ is also an interval of increase.
We can thus identify countably many, disjoint maximal intervals 
$]\alpha_i, \beta_i[$, $i\in \I^+$ of increase.
Similarly, we can identify countably many disjoint maximal intervals 
$]\gamma_i, \delta_i[$, $i\in \I^-$ of decrease.

It remains to analyze what happens at the endpoints of these intervals.
\begi
\item
If $\alpha_i\notin\S$ and, for some $\ve>0$
\bel{ep1} \mu([\alpha_i, \,\alpha_i+\ve])+\int_{\alpha_i}^{\alpha_i+\ve} {dx\over f^+(x)}~<~+\infty,\eeq
 we then consider 
the half-open interval $I^+_i\doteq [\alpha_i, \beta_i[\,$.
Otherwise, we let $I^+_i$ be an open interval:
$I^+_i\doteq \,]\alpha_i, \beta_i[\,$.
\item
If $\delta_i\notin\S$ and, for some $\ve>0$
\bel{ep2} \mu([\delta_i-\ve, \,\delta_i])-\int_{\delta_i-\ve}^{\delta_i} {dx\over f^-(x)}~<~+\infty,\eeq
 we then consider 
the half-open interval $I^-_i\doteq\, ]\gamma_i, \delta_i]\,$.
Otherwise, we let $I^-_i$ be an open interval:
$I^-_i\doteq \,]\gamma_i, \delta_i[\,$.
\endi

For each $i\in\I^+$, we now describe the increasing dynamics on the 
intervals $I^+_i$.   Given $x_0\in I_i^+$, we consider the time
$$\tau^+(x_0)~\doteq~\mu([x_0, \beta_i]) + \int_{x_0}^{\beta_i} {dy\over f^+(y)}~\in~]
0,\,+\infty].$$
We then set
\bel{S+x}
S^+_t(x_0)~\doteq~ x(t),\eeq
where $x(t)$ is implicitly defined by 
\bel{St+}\bega{cl}\ds
\mu([x_0, x(t)]) + \int_{x_0}^{x(t)} {dy\over f^+(y)}~=~t \qquad 
&\hbox{if}~t<\tau^+(x_0),\\[4mm]
x(t)~=~\beta_i\qquad&\hbox{if}~t\geq\tau^+(x_0).\enda
\eeq
The construction of the decreasing dynamics on the intervals $I_i^-$ is entirely similar.
Given $x_0\in I_i^-$, we consider the time
$$\tau^-(x_0)~\doteq~\mu([\gamma_i, x_0]) - \int_{\gamma_i}^{x_0} {dy\over f^+(y)}~\in~]
0,\,+\infty].$$
We then set
\bel{S-x}
S^-_t(x_0)~\doteq~ x(t),\eeq
where $x(t)$ is implicitly defined by 
\bel{St-}\bega{cl}\ds
\mu([x(t), x_0]) - \int^{x_0}_{x(t)} {dy\over f^-(y)}~=~t \qquad 
&\hbox{if}~t<\tau^-(x_0),\\[4mm]
x(t)~=~\gamma_i\qquad&\hbox{if}~t\geq\tau^-(x_0).\enda
\eeq
We can now combine together the above solutions, and define the semigroup
$S$ on the whole real line, as follows.
\bel{sd1}
S_t(x_0)~=~\left\{\bega{cl} x_0\qquad &\hbox{if}\quad x_0\notin\left(\bigcup_i I_i^+\right)
\cup\left(\bigcup_i I_i^-\right),\\[4mm]
S_t^+(x_0)\qquad &\hbox{if}\quad x_0\in\left(\bigcup_i I_i^+\right)\setminus
\left(\bigcup_i I_i^-\right),\\[4mm]
S_t^-(x_0)\qquad &\hbox{if}\quad x_0\in\left(\bigcup_i I_i^-\right)\setminus
\left(\bigcup_i I_i^+\right).\enda\right.\eeq
To complete the definition, it remains to define $S_t(x_0)$ in the case
$x_0\in \left(\bigcup_i I_i^-\right)\cap
\left(\bigcup_i I_i^+\right)$.   Notice that this can happen 
only if 
$$x_0~=~\alpha_i ~=~ \delta_j\,$$
where $I_i^+=[\alpha_i, \beta_i[\,$ and $I_j^-= \,]\gamma_j, \delta_j]$
are half-open intervals where the dynamics is increasing and decreasing, respectively.
By our definitions, this implies $x_0\in \Omega^*$.
Recalling {\bf (Q3)}, we thus define
\bel{sd2}
S_t(x_0)~=~\left\{ \bega{rl} S_t^+(x_0)\qquad &\hbox{if}\quad \Phi(x_0) = ~1,\\[3mm]
 S_t^-(x_0)\qquad &\hbox{if}\quad \Phi(x_0) = -1.
 \enda\right.\eeq
We claim that the above map $S$ defines a semigroup compatible with the ODE 
(\ref{ode}).    Toward this goal, we first prove

 \begin{lemma}\label{cara-so} Let $x:[0,\tau]\mapsto \R$ be a strictly increasing map defined on some interval $[0,\tau]$. Then $x(\cdot)$ is a Carath\'eodory 
solution to the Cauchy problem  (\ref{ode})-(\ref{id}) if and only if  there exists a positive, atomless  Borel measure 
$\mu$, supported on the set
$f^{-1}(0)$, such that 
\bel{incr}\int_{x_0}^{x(t)} {1\over f(s)}\, ds + \mu\bigl( [x_0, x(t)]\bigr)~=~t\eeq
%or, respectively,
%\bel{decr1}
%\int^{x_0}_{x(t)} {-1\over f(s)}\, ds + \mu\bigl( [x(t), x_0]\bigr)~=~t.\qquad\eeq
for all $t\in [0,\tau]$.
\end{lemma}
{\bf Proof.} 
{\bf 1.} Assume that $x(\cdot)$ is a Carath\'eodory 
solution to 
(\ref{ode})-(\ref{id}). Set $x_0=x(0)$ and $\Hat x=x(\tau)$.
Define  a  positive Borel  measure $\mu$ on $[x_0,\Hat x]$ by setting 
\bel{mu-t}
\mu([x_0, y])~=~\meas\bigl(\{t\in [0,t(y)]\,;~f(x(t))=0\}\bigr),
\eeq
where $t(y)$ is the unique time such that $x(t(y))=y$. Of course, the right hand side of (\ref{mu-t}) refers to Lebesgue measure.
Since $f^{-1}(0)$ is closed,  for every  $x\in ]x_0,\Hat x[\backslash f^{-1}(0)$ there exists $\delta>0$ such that $]x-\delta,x+\delta[\subset  ]x_0,\Hat x[\backslash f^{-1}(0)$ and this yields 
\[
\mu(]x-\delta,x+\delta[)~=~0.
\]
Hence $\mu$ is supported on the set $f^{-1}(0)$. By the continuity and the 
strict monotonicity property of the function $x(\cdot)$, 
the map $y\mapsto \mu([x_0,y])$ is continuous. 
This implies that $\mu$ is an atomless  measure. 

To prove (\ref{incr}) we observe that, for every $t\in [0,\tau]$, the set  
\bel{T0}
\mathcal{T}_0(t)~=~\{s\in [0,t]\,;~f(x(s))=0\}
\eeq
is~closed.
Hence $]0,t[\backslash {\mathcal T}_0(t)$  is open and can be written by as a disjoint union of countably many open intervals
\[
]0,t[\,\backslash \mathcal{T}_0(t)~=~\bigcup_{i} \,]s_i,t_i[\,.
\]
By the definitions (\ref{mu-t})-(\ref{T0}) it thus follows
\begin{eqnarray*}
t&=&\meas\Big( \{s\in [0,t]\,;~f(x(s))\not= 0\}\Big) + \meas\Big( \{s\in [0,t]\,;~f(x(s))= 0\}\Big)\\
&=&\meas\Big( ]0,t[\setminus {\mathcal T}_0(t) \Big) +\mu\bigl([x_0,x(t)]\bigr)\\
&=&\sum_{i}\int_{s_i}^{t_i}1dt+\mu\bigl([x_0,x(t)]\bigr)~=~
\sum_{i}\int_{s_i}^{t_i}{\dot{x}(t)\over f(x(t))}dt+\mu\bigl([x_0,x(t)]\bigr)\\
&=&\sum_{i}\int_{x(s_i)}^{x(t_i)}{1\over f(y)}dy+\mu\bigl([x_0,x(t)]\bigr) ~=~\int_{x_0}^{x(t)}{1\over f(y)}dy+\mu\bigl([x_0,x(t)]\bigr),
\end{eqnarray*}
proving (\ref{incr}).
\medskip

{\bf 2.} Conversely, let $\mu$ be a positive, atomless Borel measure on $[x_0, \Hat x]
=[x(0), x(\tau)]$, 
and  assume that $t\mapsto x(t)$ is a strictly increasing function, implicitly defined by (\ref{incr}) for all $t\in [0,\tau]$. For every  finite sequence of pairwise disjoint sub-intervals 
$]s_i,t_i[\,\subseteq\,]0,\tau[$, we estimate 
\begin{eqnarray*}
\sum_{i=1}^{N} |t_i-s_i|&=&\sum_{i=1}^{N}\left(\int_{x(s_i)}^{x(t_i)}{1\over f(y)}dy+\mu\bigl([x(s_i),x(t_i)]\bigr)\right)~\geq~{1\over M}\cdot \sum_{i=1}^{N}|x(t_i)-x(s_i)|.
\end{eqnarray*}
This implies that $x(\cdot)$ is absolutely continuous on $[0,\tau]$. 
In particular, $x(\cdot)$ is differentiable almost everywhere in $[0,\tau]$ and 
\bel{abs}
x(t)~=~x(s)+\int_{s}^{t}\dot{x}(r)dr\qquad\forall s<t\in [0,\tau].
\eeq
As in (\ref{Df}), we denote by $D_f$ the set of points where $f$ is discontinuous.
For any $s_0\in ]0,\tau[$ such that $x(\cdot)$ is differentiable at $s_0$ and $x(s_0)\notin f^{-1}(0)\cup D_f$, one has
\[
 \lim_{\ve\to 0+} {1\over \ve}\cdot \int^{x(s_0+\ve)}_{x(s_0)}{1\over f(y)}dy~=~{\dot{x}(s_0)\over f(x(s_0))}\,.
\]
Since $\mu$ is supported on the closed set $f^{-1}(0)$, we have $\mu([x(s_0), x(s_0+\delta)])=0$ for some $\delta>0$. Thus, (\ref{incr}) implies
\[
1~=~\lim_{\ve\to 0+}\left[{1\over \ve}\cdot \int_{x(s_0)}^{x(s_0+\ve)}{1\over f(y)}dy+ {\mu([x(s_0), x(s_0+\ve)])\over \ve}\right]~=~{\dot{x}(s_0)\over f(x(s_0))}.
\]
Since $D_f$ is countable, one has 
\bel{x's}
\dot{x}(s)~=~f(x(s))~~~ \hbox{for a.e.}~~s\in [0,\tau]\backslash \mathcal{T}_0(\tau).
\eeq
Given any $t\in ]0,\tau[$, the open set $]x(0),x(t)[\backslash f^{-1}(0)$ can be written as a disjoint union of countably many open intervals $]\alpha_k,\beta_k[$. Setting $s_k=x^{-1}(\alpha_k)$ and $t_k=x^{-1}(\beta_k)$, we obtain
\[
\bigcup_{k\geq 1}]s_k,t_k[~=~]0,t[\,\backslash \mathcal{T}_0(t).
\]
Thus, (\ref{x's}) and (\ref{abs}) yield
\begin{eqnarray*}
x(t)-x(0)&=&\sum_{k\geq 1} (x(t_k)-x(s_k))~=~\sum_{k\geq 1}\int_{s_k}^{t_k}\dot{x}(s)ds\\
&=&\int_{\bigcup_{k\geq 1} ]s_k,t_k[}f(x(s))ds~=~\int_{]0,t[\backslash \mathcal{T}_0(t)}f(x(s))ds~=~\int_{0}^{t}f(x(s))ds
\end{eqnarray*}
showing that $x(\cdot)$ is a Carath\'eodory solution. This completes the proof.
\endproof
\medskip

We are now ready to state the main result of this section.

\begin{theorem}\label{car_sem}  Let $f:\R\mapsto\R$ satisfy the assumptions {\bf (A1)-(A2)}. The following statements are equivalent.
\begi
\item[(i)] The map $S:\R\times \R_+\to\R$ is a deterministic semigroup compatible with the ODE (\ref{ode}).
\item [(ii)] There exist a positive atomless measure $\mu$ supported on the 
set $f^{-1}(0)$, a countable set of points $\S\subseteq f^{-1}(0)$, and a map $\Phi:\Omega^*\mapsto\{-1,1\}$
as in {\bf (Q1)--(Q3)} such that $S$ coincides with the corresponding 
semigroup constructed at (\ref{sd1})-(\ref{sd2}).
\endi
\end{theorem}

{\bf Proof of the implication (ii)$\implies$(i).} Given 
the measure $\mu$, the set $\S$ and the map $\Phi$, 
we will show that 
the map $S:\R\times\R_+\mapsto\R$ constructed in (\ref{sd1})-(\ref{sd2})
is indeed a semigroup compatible with the ODE (\ref{ode}).
\medskip

{\bf 1.}
Let $x_0\in \R$ be given.  We first show that the map $t\mapsto x(t)= S_t(x_0)$ provides
 a solution to the Cauchy problem (\ref{ode})-(\ref{id}). 
Various cases must be considered.

{\bf CASE 1.}  If $x_0\notin\left(\bigcup_i I_i^+\right)
\cup\left(\bigcup_i I_i^-\right)$ then $S_t(x_0)=x_0$. In this case we need
to show that $f(x_0)=0$. If, on the contrary, $f(x_0)\not= 0$, then by {\bf (A2)} we would have either 
\bel{ctr1}0~<~\delta~\leq~ f(x)~\leq~ M\eeq for all 
$x\in \cN_\ve\doteq \,] x_0, x_0+\ve]$,
or    else \bel{ctr3}-M~\leq ~f(x)\leq ~-\delta~<0 \eeq for all 
$x\in  \cN_\ve \doteq \, [x_0-\ve, x_0[$, for a suitably small $\ve>0$.

We now recall that the measure $\mu$ and the countable set 
$\S$ in {\bf (Q1)-(Q2)} satisfy
$$Supp(\mu)\cup \S~\subseteq~ f^{-1}(0).$$
Since the set $f^{-1}(0)$ is closed and 
does not contain $x_0$, by possibly choosing a  smaller $\ve>0$ we can assume that 
\bel{ctr2}\S\cap\cN_\ve~=~\emptyset,\qquad \mu(\cN_\ve)=0.\eeq
Together, (\ref{ctr1})--(\ref{ctr2}) imply that either a strictly increasing trajectory
or a strictly decreasing trajectory 
starts from $x_0$.  Hence $x_0\in\left(\bigcup_i I_i^+\right)
\cup\left(\bigcup_i I_i^-\right)$, against the assumption.

{\bf CASE 2.}   Next, assume that $x_0\in I_i^+\setminus
\left(\bigcup_j I_j^-\right)$, for some $i\in\I^+$ and $I_i^+= [\alpha_i, \beta_i[$.  
Then the trajectory $t\mapsto x(t)= S_t^+(x_0)$ is implicitly defined by (\ref{incr}).
By Lemma~\ref{cara-so}, we already know that $x(\cdot)$ is a 
Carath\'eodory solution to the ODE (\ref{ode}) for $t\in [0, \tau(x_0)]$.

It remains to show that, in the case $\tau(x_0)<+\infty$, the function $x(t)=\beta_i$
is still a solution of (\ref{ode}) for all $t\in [\tau(x_0), +\infty[\,$.
Of course, this will be the case iff $f(\beta_i)=0$. 

Assume that, on the contrary, $f(\beta_i)\not=0$.  Since $f^{-1}(0)$ is closed, 
there exists a neighborhood $\cN$ of $\beta_i$ such that (\ref{ctr2}) holds.
Moreover, since $f(\beta_i-)\geq 0$, by {\bf (A2)} it follows $f(\beta_i+)>0$.
This implies that, for some $\ve>0$ small enough, the inequalities (\ref{ctr1}) hold 
for all $x\in \cN_\ve\doteq \,] \beta_i, \beta_i+\ve]$. Moreover, (\ref{ctr2}) holds as well.
As a consequence, there exists a strictly increasing trajectory starting at $\beta_i$.
This yields a contradiction with the maximality of the interval $I_i^+$.
\v
{\bf CASE 3.}  If $x_0\in I_i^-\setminus
\left(\bigcup_j I_j^+\right)$, for some $i\in\I^-$, 
then the same arguments used in Case 2 show that $t\mapsto 
x(t) = S^-_t(x_0)$ is a Carath\'eodory solution.
\v
{\bf CASE 4.}  If $x_0\in \left(\bigcup_i I_i^-\right)\cap
\left(\bigcup_i I_i^+\right)$, then by (\ref{sd2}) we 
are back to Case 2 or Case 3.
\v
{\bf 2.} It remains to prove that the semigroup property $S_s(S_t(x_0)) = S_{s+t}(x_0)$
is satisfied.

Referring to the above four cases, in Case 1  this is trivial, because $S_t(x_0)=x_0$
for all $t\geq 0$. 

In Case 2, as long as $t+s\leq \tau^+(x_0)$, the semigroup property is a straightforward consequence of the identity
$$\bega{rl} t+s&\ds =~\mu([x_0, x(t+s)]) + \int_{x_0}^{x(t+s)} {dy\over f^+(y)}\\[4mm]
&\ds=~\mu([x_0, x(t)]) + \mu([x(t), x(t+s)])+ \int_{x_0}^{x(t)} {dy\over f^+(y)}+\int_{x(t)}^{x(t+s)} {dy\over f^+(y)}.\enda$$
On the other hand, for $0<t< \tau(x_0)\leq t+s$, the identity
$$\tau^+(S_t(x_0))~=~\tau^+(x_0)-t$$
implies 
\bel{sbi}S_{t+s} (x_0) ~=~\beta_i~=~S_s( S_t(x_0)).\eeq
Finally, if $t\geq\tau^+(x_0)$, then again (\ref{sbi}) trivially holds.

The remaining Cases 3 and 4 are entirely similar.
This concludes the proof of the implication (ii)$\implies$(i).
\v
{\bf Proof of the implication (i)$\implies$(ii).} Assuming that $S:\R\times\R_+\mapsto\R$ is a deterministic semigroup compatible with the ODE (\ref{ode}),
in the remaining part of the proof 
we shall construct a measure $\mu$, a countable set $\S$, and a map 
$\Phi$ so that $S$ coincides with the corresponding 
semigroup defined at (\ref{sd1})-(\ref{sd2}).  
\medskip

{\bf 3.}
 In analogy with Definition~\ref{d:increase},
we introduce:

\begin{definition}\label{d:incdyn} Given the semigroup $S:\R\times\R_+\mapsto\R$, we say that an  interval $J= [a,b[$ or $J=\,]a,b[$ (open to the right) is a 
{\bf domain of increasing dynamics}
if, for every $x_0, \Hat x\in J$ with $x_0<\Hat x$ there exists $t>0$ such that 
$\Hat x ~=~S_t x_0$.

Similarly, we say that an  interval $J= ]a,b]$ or $J=\,]a,b[$ (open to the left) is a 
{\bf domain of decreasing dynamics}
if, for every $x_0, \Hat x\in J$ with $\Hat x < x_0$ there exists $t>0$ such that 
$\Hat x ~=~S_t x_0$.
\end{definition}
We observe that, if two  domains of increasing dynamics $J, J'$ have nonempty
intersection, then the union $J\cup J'$ is also an interval of increasing dynamics.
We can thus partition the real line as
\bel{lp}
\R~=~\Omega_S^+\cup\Omega_S^-\cup\Omega_S^0\eeq
where
\begi
\item $\Omega_S^+=\bigcup_{j\in \J^+}  J_j^+$ is the  union of countably many 
disjoint intervals $J_j^+= [a_j, b_j[$ or $J_j^+= \,]a_j, b_j[$ open to the right, where the dynamics is increasing and  $b_j$ can not be crossed from the left, i.e. $S_t(b_j-\ve)\leq b_j$ for all $t\geq 0$, $\ve\in ]0,b_j-a_j[$.
\item $\Omega_S^-=\bigcup_{j\in \J^-}  J_j^-$ is the  union of countably many 
disjoint intervals $J_j^-=\,]c_j, d_j]$ or $J_j^-=\,]c_j, d_j[\,$ open to the left, where the dynamics is decreasing and  $c_j$ can not be crossed from the right, i.e. $S_t(c_j+\ve)\geq c_j$ for all $t\geq 0$, $\ve\in ]0,d_j-c_j[$.

\item $\Omega_S^0\subseteq f^{-1}(0)$ is the set of points $x_0$ such that 
$S_t(x_0)=x_0$ for all $t>0$.
\endi
{\bf 4.}  We shall construct the measure $\mu$ separately on each interval
$J_j^+$, $J_j^-$.   As in the proof of Lemma~\ref{cara-so}, we define a  positive,  Borel  measure on $J_j^+$ by setting 
for every compact subinterval
$[x_0, \Hat x]\subset J_j^+$ 
\bel{mud3} \mu\bigl([x_0, \Hat x]\bigr)~\doteq~\meas
\Big( \bigl\{ t\in [0,\tau]\,;~f(S_t(x_0))=0\bigr\}\Big)\eeq
where $\tau>0$ is such that  $\Hat x= S_\tau (x_0)$. It is clear that $\mu$ is supported on $f^{-1}(0)$. By the semigroup property of $S$, the map $t\mapsto S_t(x_0)$ is strictly increasing. Thus, as in the proof of Lemma~\ref{cara-so}, the measure $\mu$ is also atomless. Similarly, one can define the measure $\mu$ on $J_j^-$.

Next, we define the countable set
\bel{Sdef3}
\S~\doteq~\Big(\bigl\{ a_j\,,\,b_j:j\in \J^+\} \cup \bigl\{ c_j\,,\,d_j: j\in \J^-\}
\Big) \cap~\Omega_S^0
\eeq
containing all endpoints of all intervals $J_j^+, J_j^-$  whose trajectory 
remains constant.

Finally,  the discrete function $\Phi:\Omega^*\mapsto \{-1,1\}$ is defined as follows.
Let $x_0\in \Omega^*$. Three cases are considered.
\begi
\item If $x_0= a_j$ for some half-open interval $J_j^+=[a_j, b_j[\,$, then 
we set
$\Phi(x_0)=1$.
\item  If $x_0= d_j$ for some half-open interval $J_j^-=\,]c_j, d_j]\,$, then 
we set
$\Phi(x_0)=-1$.
\item In all other cases, the choice of $\Phi(x_0)$ is immaterial, because it has
no effect on the definition of the semigroup.
To fix the ideas, we may thus set $\Phi(x_0)=1$.
\endi

Recalling Definition \ref{d:increase}, let $I_i^+$, $i\in \I^+$, and $I_i^-$, $i\in \I^-$, be respectively the domains of increase and the domains of decrease corresponding to $\mu$, $\S$ and $\Phi$. To complete the proof, we first show  that 
\bel{conci}
\{]\alpha_i,\beta_i[:i\in\I^+\}~=~\{]a_j,b_j[:j\in\J^+\}\quad\mathrm{and}\quad \{]\gamma_i,\delta_i[:i\in\I^-\}~=~\{]c_j,d_j[:j\in\J^-\}
\eeq
For any $j\in \J^+$, by  Lemma~\ref{cara-so} and (\ref{mud3})-(\ref{Sdef3}), one has that $]a_j,b_j[$ satisfies (\ref{din}) and this implies that $]a_j,b_j[\subset I^+_{i}$ for some $i\in\mathcal{I}^+$. Moreover, if $a_j\in J_j^+$ then $a_j\notin \Omega_S^0$ and 
\[
\mu([a_j,a_j+\ve[)+\int_{a_j}^{a_j+\ve}{1\over f^+(y)}dy~=~t_{\ve}~<~\infty\qquad\mathrm{with}\qquad S_{t_{\ve}}(a_j)~=~a_j+\ve
\]
for $0<\ve<b_j-a_j$. Hence, $a_j\in I_i^+$ and this yields $J_j^+\subseteq I_i^+$. Furthermore, we prove that $b_j=\beta_i, a_j=\alpha_i$ and this implies  
\[
]\alpha_i,\beta_i[~=~]a_j,b_j[.
\]
If  $b_j<\beta_i$ then $b_j$ is reachable from the left and  the semigroup property of $S$ implies that $b_j\in \Omega_S^0$. In particular, $b_j\in\S\cap I_i^+$ and this yields a contradiction. Now assume that $\alpha_i<a_j$. Two cases are considered:
\begin{itemize}
\item If $a_j\notin J_{j}^+$ then $a_j\in \Omega^0_S$ since $\Omega^-_S\cap I_i^+=\emptyset$. Thus, $a_j\in \S\cap I_i^+$ and this yields a contradiction. 
\item Otherwise, $a_j\in J_{j}^+$ can not be reached from the left. Thus, either $\mu([a_j-\ve,a_j])=+\infty$ for some $\ve\in ]0,a_j-\alpha_i[$ or there exists $b_{j'}\in \Omega_S^0\cap I_i^+$. Both cases yield a contradiction. 
\end{itemize}
Vice-versa, for any $i\in \I^+$, by the definition \ref{d:incdyn} and Lemma \ref{cara-so}, one has that
$$
\left(\Omega_S^-\bigcup\mathrm{int}(\Omega_S^0)\right)\cap~ ]\alpha_i,\beta_i[~=~\emptyset.
$$
In particular,  $]a_j,b_j[\cap ]\alpha_i,\beta_i[\neq \emptyset$ for some $j\in \J^+$. Since $]a_j,b_j[\subseteq ]\alpha_{i'},\beta_{i'}[$ for some $i'\in \I^+$ and $\{I_{k}\}_{k\in\I^+}$ are disjoint, we have that $]a_j,b_j[=]\alpha_i,\beta_i[$, and the first equality of (\ref{conci}) holds. Similarly, one can show that the second equality of (\ref{conci}). Finally, we need to show that 
\begin{itemize}
\item [(i)] If $\alpha_i \in I_{i}^+\backslash J_j^+$ and $]\alpha_i,\beta_i[=]a_j,b_j[$ for some $j\in\J^+$ then 
$$
 \alpha_i~=~d_k~\in~\Omega^*\cap J_k^-\qquad\mathrm{for~some}~k\in \J^-.
$$ 
\item [(ii)]  If $\delta_i \in I_{i}^-\backslash J_j^-$ and $]\gamma_i,\delta_i[=]c_j,d_j[$ for some $j\in\J^-$ then 
$$
 \delta_i~=~a_k~\in~\Omega^*\cap J_k^+\qquad\mathrm{for~some}~k\in \J^+.
$$ 
\end{itemize}
Since (i) and (ii) are entirely similar, let us prove (i). Observe that $\alpha_i=a_j$ can not be in $\Omega^0_S$, it holds that $\alpha_i=d_{k}$ for some $k\in \J^-$. On the other hand, one can follow the above argument to show that $J^-_{k}\subseteq I^-_{\ell}$ for some $\ell\in\I^-$.  Thus, $\alpha_i=d_k\in I^+_i\cap I^-_{\ell}\subseteq\Omega^*$ and this yields (i).
\endproof

\section{A Markov  process whose sample paths solve the ODE}
\setcounter{equation}{0}
\label{s:5}
In this section, we initiate the study of Markov semigroups whose sample paths
satisfy the ODE (\ref{ode}) with probability one.
Our eventual goal is to show that each of these random 
semigroups can be obtained starting with a deterministic semigroup and performing the two modifications described at {\bf (Q4)-(Q5)}
in the Introduction.   
%Namely,
%\begi
%\item  For a countable set $\S^*$ of zeroes of $f$, we assign an exponential
%waiting time.
%\item  For the set of isolated points $\Omega^*$ in (\ref{O*}), we assign the probabilities of going left or right.
%\endi
In this section, we show that  a Markov process can be uniquely determined
by the ODE (\ref{ode}) and the assignments {\bf (Q1)-(Q2)} and {\bf (Q4)-(Q5)},
namely:
\begi
\item A positive, atomless Radon measure supported on $f^{-1}(0)$.
\item A countable set $\S\subseteq f^{-1}(0)$ where the dynamics stops.
\item A countable set $\S^*\subseteq f^{-1}(0)$ and, for each $y_j\in \S^*$, a 
number $\lambda_j>0$ characterizing
the random  exponential 
waiting time, when the trajectories reach the point $y_j\in \S^*$.
\item A map $\Theta:\Omega^*\mapsto [0,1]$ determining the probability of moving 
to the right or to the left, from an initial point $z_k\in \Omega^*$.
\endi

We start by selecting an underlying probability space ${\mathcal W}$ such that,
as $\omega\in {\mathcal W}$, the countably many random variables $Y_j(\omega)\in \R_+$ and $Z_k(\omega)\in \{-1,1\}$ are independent, with distributions 
\bel{pxj}
\P\{ Y_j>s\}~=~e^{-\lambda_j s}\quad\forall s>0,\quad
\qquad \quad \P\{ Z_k=1\} ~=~\Theta(z_k).\eeq

Next, we recall that, given the measure $\mu$ and the countable set $\S$,
as in Section~\ref{s:3} one can uniquely determine the sets 
\bel{om+-}\Omega^+~=~\bigcup_{i\in\I^+} I_i^+\,,\qquad\qquad 
\Omega^-~=~\bigcup_{i\in\I^-} I_i^-\eeq
consisting of countable unions of disjoint intervals where the dynamics can increase or decrease, respectively.
For $x_0\in \Omega^+$, a trajectory $t\mapsto S_t^+(x_0)$ was defined at 
(\ref{S+x})-(\ref{St+}), while for 
$x_0\in \Omega^-$, a trajectory $t\mapsto S_t^-(x_0)$ was defined at 
(\ref{S-x})-(\ref{St-}).   

To construct a Markov process with sample paths $t\mapsto X(t,\omega)$
satisfying (\ref{ode}), 
for  every
$t>0$ we need to define the transition probabilities
$$P_t(x_0, A)~=~\P\Big\{X(t,\omega)\in A\,\Big|~~X(0,\omega)= x_0\Big\},$$
for any initial point $x_0\in \R$ and any Borel set $A\subset\R$.
Recalling our construction of a deterministic flow 
at (\ref{sd1})-(\ref{sd2}), this can be done as follows.
\begi
\item[(i)] If $x_0\notin\Omega^+\cup\Omega^-$, then 
$X(t, x_0, \omega) = x_0$ for every $t\geq 0$ and $\omega\in {\mathcal W}$.
Hence 
\bel{TP1} P_t(x_0, \{x_0\})~ = ~1\qquad \forall ~t\geq 0.\eeq
In other words, all  trajectories starting at $x_0$ remain constant.
\item[(ii)] If $x_0\in  \Omega^+\setminus\Omega^-$, a random trajectory starting at 
$x_0$ will have the form
\bel{X+}
X(t, x_0, \omega)~=~S^+_{T^+(t,x_0,\omega)} (x_0)\,,\eeq
where the time $T^+$ along the trajectory is a random variable with distribution 
\bel{T+} {\P\bigl\{ T^+(t, x_0, \omega)<s\bigr\}}~=~\P\left\{ s+
\sum_{y_j\in S^*\cap [ x_0, \, S^+_s(x_0)]}  Y_j (\omega)~>~t\right\}.\eeq
Note that (\ref{T+}) accounts for the (possibly countably many) waiting times
when the trajectory crosses one of the points $y_j\in \S^*$.
The transition probabilities are thus given by
\bel{TP2}\bega{l}
P_t\bigl( x_0,\, \,]-\infty, x_0[\bigr) ~=~0,\\[3mm]
\ds P_t\bigl( x_0,\, [x_0, S_s^+(x_0)]\bigr)~=~\P\left\{ s+\sum_{y_j\in S^*\cap [ x_0, \, S^+_s(x_0)]}  Y_j (\omega)~\geq ~t\right\}.
\enda \eeq

\item[(iii)] Similarly, for an initial state  $x_0\in  \Omega^-\setminus\Omega^+$, 
a random trajectory starting at 
$x_0$ will have the form
\bel{X-}
X(t, x_0, \omega)~=~S^-_{T^-(t,x_0,\omega)} (x_0)\,,\eeq
where the time $T^-$ along the trajectory is a random variable with distribution 
\bel{T-} {\P\bigl\{ T^-(t, x_0, \omega)<s\bigr\}}~=~\P\left\{ s+
\sum_{y_j\in S^*\cap [ S^-_s(x_0),\, x_0]}  Y_j (\omega)~>~t\right\}.\eeq
The transition probabilities are thus given by 
\bel{TP3}\bega{l}
\ds
P_t\bigl( x_0,\, \,]x_0, +\infty[\bigr) ~=~0,\\[3mm]
\ds P_t\bigl( x_0,\, [S_s^-(x_0), x_0]\bigr)~=~\P\left\{ s+\sum_{y_j\in S^*\cap 
[ S^-_s(x_0), \, x_0]}  Y_j (\omega)~\geq ~t\right\}.\enda\eeq
%\bel{X+}
%X(t, x_0, \omega)~=~S^-_{T(t,x_0,\omega)} (x_0)\,,\eeq
%where the time $T^-$ along the trajectory  is a random variable with distribution 
%\bel{T-} \P\bigl\{ T^-(t, x_0, \omega)>s\bigr\}~=~\P\left\{ s+
%\sum_{y_j\in \S^*\cap [ \, S^-_s(x_0), x_0]}  Y_j ~>~t\right\}.\eeq
\item[(iv)]
To complete the definition, it remains to define the transition probabilities in the case
$x_0 \in\Omega^+\cap\Omega^-\subseteq\Omega^*$.     In this case,
by construction we have $x_0= z_k$ for some $k$.
We then define
the random variable 
$X(t, x_0,\omega)$ by setting
\bel{X**} 
X(t, x_0,\omega)~=~\left\{ \bega{rl} S^+_{T^+(t, x_0, \omega)}(x_0)\quad
&\hbox{if}\quad Z_k(\omega)=1\,,\\[4mm]
S^-_{T^-(t, x_0, \omega)}(x_0)\quad
&\hbox{if}\quad Z_k(\omega)= -1\,.\enda\right.\eeq
By (\ref{pxj}), its distribution satisfies
\bel{TP4}P_t(x_0,A)~=~
\Theta(x_0)\cdot 
\P\bigl\{S^+_{T^+(t, x_0, \omega)}\in A\bigr\} + (1-\Theta(x_0))\cdot 
\P\bigl\{S^-_{T^-(t, x_0, \omega)}\in A\bigr\} \eeq
for every Borel set $A\subseteq\R$.
\endi

From the above construction it is clear that, in all cases (i)--(iv), all sample paths
 $t\mapsto X(t, x_0,\omega)$ are Carath\'eodory solutions to the  Cauchy
 problem  (\ref{ode})-(\ref{id}). Indeed, we are only adding a countable number of waiting times, when the random trajectories reach one of the points $y_j\in \S^*\subseteq 
 f^{-1}(0)$. 

The next result shows that the above  
transition probabilities define a Markov semigroup.
\v
\begin{proposition}\label{p:5} The transition probability kernels 
$(P_t)_{t\geq 0}$ in (\ref{TP1}), (\ref{TP2}), (\ref{TP3}), (\ref{TP4}) define a
continuous  Markov process.
%, whose sample paths are Carath\'eodory solutions to
%the ODE (\ref{ode}).    we have $P_0(x,\cdot)=\delta_x$ and 
Indeed,  they satisfy the Chapman-Kolmogorov equation
\begin{equation}\label{ckeq}
\int P_t(z,A) P_s(x_0,dz)~=~P_{s+t}(x_0,A)
\end{equation}
for all $ s,t>0$, $ x_0\in\R$ and every Borel set $ A\subseteq \R$.
\end{proposition}

{\bf Proof.} 
{\bf 1.} In Case (i) the conclusion is trivially true. 
We thus focus on Case (ii). Suppose that $x_0\in I^+_{i}$ for some $i\in \I^+$.
Observe that
$X(s,x_0,\omega)=z\doteq  S_\eta^+(x_0)$ if and only if
\bel{u1} \eta+\sum_{y_j\in S^*\cap [ x_0, \, S^+_\eta(x_0)[}  
Y_j (\omega)~\leq ~s ~\leq~ \eta+\sum_{y_j\in S^*\cap [ x_0, \, S^+_\eta(x_0)]}  
Y_j (\omega)\eeq
This leads to two different expressions, 
depending on whether $z\in S^*$ or $z\notin S^*$.
We claim that, for any Borel set $A\subset\R$,
\bel{u2}
P_{t+s}(x_0, A)~\ds=~
\sum_{y_j\in \S^*\cap [x_0, S^+_s(x_0)]} P_s\bigl(x_0, \{y_j\}\bigr)\cdot P_t(y_j, A)
+\int_{z\in [x_0, S^+_s(x_0)] \setminus \S^*}P_s(x_0, dz)P_t(z, A).
\eeq
 It is sufficient to prove (\ref{u2}) assuming 
 $x_0\in I_i^+$ and   $A=[x_0,b]\subseteq I^+_i$, where $I_i^+$ is one of the intervals in (\ref{om+-})
where the dynamics is increasing.
Consider the times $\eta_j, \eta_z\geq 0$ such that
$y_j= S^+_{\eta_j}(x_0)$ and $z= S^+_{\eta_z}(x_0)$. 
\begi
\item The time when the random trajectory $t\mapsto X(t, x_0,\omega)$ crosses $z$ is
given by
\[
T^+_z(\omega)~=~\eta_z+\sum_{y_j\in S^*\cap [x_0,z]} Y_j(\omega).
%\qquad \quad C~=~\left\{\omega:t+s\leq T^+_b(\omega)\right\}.
\]
\item The times when the random trajectory  $t\mapsto X(t, x_0,\omega)$ reaches $y_j$,
and then leaves $y_j$, are given by
\[
T^-_{j}(\omega)~=~\eta_j+\sum_{y_{k}\in S^*\cap [ x_0, y_j[}  Y_{k} (\omega),\qquad T^+_{j}(\omega)~=~\eta_j+\sum_{y_{k}\in S^*\cap [ x_0, y_j]}  Y_{k} (\omega).
\]
\endi
In view of the above definition, the left hand side of (\ref{u2}) can be written as 
\begin{multline}\label{u2'}
\ds P_{t+s}\bigl(x_0, [x_0,b]\bigr)~=~\P\Big\{ T_b^+(\omega)\geq t+s\Big\}\\
 =~\P\Big\{ T_b^+(\omega)\geq t+s~~ \hbox{and}~~ X(s,x_0,\omega)\in S^*\Big\}+ \P\Big\{ T_b^+(\omega)\geq t+s~~ \hbox{and}~~ X(s,x_0,\omega)\notin S^*\Big\}\\
 =~\sum_{y_j\in S^*\cap [x_0,S^+_s(x_0)]} \P\Big\{ T_b^+(\omega)\geq t+s~~ \hbox{and}~~ 
X(s,x_0,\omega)= y_j \Big\} \\+ \P\Big\{ T_b^+(\omega)\geq t+s~~ \hbox{and}~~ X(s,x_0,\omega)\notin S^*\Big\}~=~I+II.
\end{multline}
%\bel{u2'}\bega{l}
%\ds P_{t+s}\bigl(x_0, [x_0,b]\bigr)~=~\P\Big\{ T_b^+(\omega)\geq t+s\Big\}\\[4mm]
%\quad =~ \P\Big\{ T_b^+(\omega)\geq t+s~~ \hbox{and}~~ X(s,x_0,\omega)\in S^*\Big\}+ \P\Big\{ T_b^+(\omega)\geq t+s~~ \hbox{and}~~ X(s,x_0,\omega)\notin S^*\Big\}.\enda\eeq
%\sum_{y_j\in S^*\cap [x_0,S^+_s(x_0)]}\P\left\{C\bigcap\big\{s\in [T^+_{j-},T^+_{j^+}]\big\}\right\}+\P\left\{C\cap \left\{S^+_{T^+(s,x_0,\omega)}(x_0)\notin S^* \right\}\right\}.
%\eeq
\v
{\bf 2.} 
For any $y_j\in S^*\cap [x_0,S^+_s(x_0)]$, we have 
\[
\bega{l}\ds\P\Big\{ T_b^+(\omega)\geq t+s~~ \hbox{and}~~ X(s,x_0,\omega)= y_j\Big\}\\[4mm]
\ds \qquad =~\P\Big\{ T_b^+(\omega)\geq t+s~~ \hbox{and}~~ 
T^-_{j}(\omega)\leq s\leq T^+_{j}(\omega) \Big\}\\[4mm]
\ds\qquad =~\int_{0}^s{d\over dr_1} \P\Big\{T^-_{j}(\omega)\in [0, r_1]\Big\}
\int_{s-r_1}^{\infty}{d\over dr_2}\P\Big\{Y_j(\omega)\in [0, r_2]\Big\}\\[3mm]
\qquad\qquad \qquad\qquad\qquad 
\ds \cdot \P\Big\{t+s-r_1-r_2\leq T^+_b(\omega)-T^+_{j}(\omega)\Big\}
\, dr_2 \,dr_1\\[4mm]
\ds \qquad  =~\int_{0}^s{d\over dr_1} \P\Big\{T^-_{j}(\omega)\in [0, r_1]\Big\}\\[3mm]
\qquad\qquad \qquad\qquad\qquad 
\ds \cdot 
 \int_{s-r_1}^{\infty}\lambda_{j}e^{-\lambda_{j}r_2}\cdot \P\left\{t+s-r_1-r_2\leq T^+_b(\omega)-T^+_{j}(\omega)\right\}dr_2\, dr_1\,.
\enda\]
Performing the change of variable $\tau= r_2-(s-r_1)$, one obtains
\[\bega{l}
\ds \P\Big\{ T_b^+(\omega)\geq t+s~~ \hbox{and}~~ X(s,x_0,\omega)= y_j\Big\}
\\[4mm]
\ds =~\int_{0}^s{d\over dr_1} \P\Big\{T^-_{j}(\omega)\in [0, r_1]\Big\}\int_{0}^\infty\lambda_{j}e^{-\lambda_{j}(\tau+s-r_1\}}\cdot \P\left\{t-\tau\leq T^+_{b}-T^+_{j}\right\}d\tau\, dr_1\\[4mm]
\ds =~\left(\int_{0}^s{d\over dr_1} \P\Big\{T^-_{j}(\omega)\in [0, r_1]\Big\}\cdot e^{-\lambda_j(s-r_1)}\, dr_1\right)\\[4mm]
\ds \qquad\qquad\qquad\qquad   \cdot 
\left(\int_{0}^\infty\lambda_{j}e^{-\lambda_{j}\tau}\cdot \P\left\{t-\tau\leq 
T_b^+(\omega)-T_j^+(\omega)\right\}d\tau\right)
\\[4mm]
\ds =~\int_{0}^s{d\over dr_1} \P\Big\{T^-_{j}(\omega)\in [0, r_1]\Big\}\cdot
\P\Big\{Y_j(\omega)\geq s-r_1\Big\}\, dr_1\\[3mm]
\ds\qquad\qquad\qquad\qquad  \cdot 
\int_{0}^\infty{d\over d\tau} \P\Big\{Y_j(\omega)\in [0,\tau]\Big\} 
\cdot \P\Big\{t-\tau\leq T^+_{b}(\omega)-T^+_{j}(\omega)\Big\}\, d\tau \\[4mm]
\ds =~\P\Big\{ s\in \bigl[T^-_j(\omega),T^+_{j}(\omega)\bigr]\Big\}\cdot 
\P\Big\{ t\leq T^+_b(\omega)-T^-_j(\omega)\Big\}\\[4mm]
=~P_s\bigl(x_0,\{y_{j}\}\bigr)\cdot 
P_t\bigl(y_{j}, [x_0,b]
\bigr).
\enda \]
Summing over all $y_j\in S^*$ we thus obtain
\bel{I}
I~=~\sum_{y_j\in S^*\cap [x_0,S^+_s(x_0)]}P_s(x_0,y_{j})\cdot P_t\bigl(y_{j}, [x_0,b]\bigr).
\eeq

{\bf 3.}
Next, for a fixed integer $N\geq 1$,  we consider the following partition of $[x_0,b]$ 
\bel{Zi}
z_0~\doteq~x_0,\qquad z_i~\doteq~\sup\Big\{z\in[x_0,b]:P_t\bigl(z,[x_0,b]\bigr)
\geq 1-{i\over N}\Big\},\eeq
and define the open intervals
$$ J_i~\doteq~\bigl]z_{i-1},z_i\bigr[
\,,\qquad\qquad i\in \{1,\ldots, N\}.$$ 
Notice that some of these intervals $J_i$ may be empty. 
For each $i\in \{1,\dots, N\}$, we estimate
$$
\bega{l}\ds
\P\Big\{t+s\leq T^+_b(\omega)~\Big|~X(s,x_0,\omega)\in J_i\setminus 
S^*\Big\}\\[4mm]
\qquad \ds \leq~\sup_{z\in J_i\smallsetminus S^*}\P\Big\{t+s\leq T^+_b(\omega)~\Big|~X(s,x_0,\omega)=z\Big\}\\[4mm]
\ds\qquad =~\sup_{z\in J_i\setminus S^*}\P\Big\{t+s\leq T^+_b(\omega)~\Big|~s= T^+_{z}(\omega)\Big\}~=~\sup_{z\in J_i\setminus S^*}\P\Big\{ T^+_{b}(\omega)-T^+_{z}(\omega)
\geq t\Big\}\\[4mm]
\qquad\ds =~\sup_{z\in J_i\setminus S^*}P_t\bigl(z, [x_0,b]\bigr)~
\leq~\inf_{z\in J_i\setminus S^*}P_t\bigl(z, [x_0,b]\bigr)+{1\over N}.
\enda
$$
An entirely similar argument yields
\[
\P\Big\{t+s\leq T^+_b(\omega)~\Big|~X(s,x_0,\omega)\in J_i\setminus 
S^*\Big\}~\geq~\sup_{z\in J_i\setminus S^*}P_t(z, [x_0,b])-{1\over N}\,.
\]
Summing over $i=\{1,\ldots,N\}$, one obtains
$$\bega{l}\ds
\P\Big\{t+s\leq T^+_b(\omega) ~~\hbox{and}~~X(s, x_0, \omega)\notin S^* \Big\}\\[4mm]
\qquad\qquad\qquad \ds = ~\sum_{i=1}^N\P\Big\{t+s\leq T^+_b(\omega) ~~\hbox{and}~~X(s, x_0, \omega)\in J_i\smallsetminus S^*\Big\}\\[4mm]
\qquad\qquad\qquad\ds=~\sum_{i=1}^N \P\Big\{t+s\leq T^+_b(\omega)~\Big|~X(s,x_0,\omega)\in J_i\setminus S^*\Big\}\cdot P_s\bigl(x_0,J_i\backslash S^*\bigr)\\[4mm]
\qquad\qquad\qquad\ds \in~\sum_{i=1}^{N}\int_{z\in J_i\backslash S^*}P_s(x_0,dz)P_t\bigl(z,[x_0,b]
\bigr)+\Big[-{1\over N},{1\over N}\Big]. 
\enda $$
By letting $N$ to $+\infty$, we thus obtain
\bel{II}
II~=~\int_{z\in [x_0, S^+_s(x_0)] \setminus \S^*} P_t(z, [x_0,b]) \cdot P_s(x_0, dz).
\eeq
Together with (\ref{I}) and (\ref{u2'}), this yields (\ref{u2}).

Case (iii) is proved in an analogous way. 
\v

{\bf 4.} In Case (iv), we show that 
the result follows from the same analysis as  in Cases (ii) and (iii). 
Indeed, it will be sufficient to prove (\ref{ckeq}) for $A=[x_0,b]$, while the case $A=[c,x_0]$ is entirely similar.

 If $x_0\notin\S^*$ then there is no random waiting time at $x_0$, and all trajectories
 immediately start moving (either to the left or to the right of $x_0$).
The same argument used in Case (ii) now yields
\[\bega{l}\ds
\P\left\{S^+_{T^+(t+s,x_0,\omega)}\in [x_0,b]\right\}\\[4mm]
\ds\qquad =~\int_{[x_0,b]}\P\left\{S^+_{T^+(t,z,\omega)}\in [x_0,b]\right\}\cdot {d\over dz}\P\left\{S^+_{T^+(s,x_0,\omega)}\in [x_0,z]\right\}dz.\enda
\] 
From (\ref{TP4}) it thus  follows
$$\bega{l}
\ds\int P_t\bigl(z,[x_0,b]\bigr) P_s(x_0,dz)~=~
\int_{[x_0,b]} P_t\bigl(z,[x_0,b]\bigr) P_s(x_0,dz)\\[4mm]
\qquad \ds=~\Theta(x_0) \cdot\int_{[x_0,b]} 
\P\left\{S^+_{T^+(t,z,\omega)}\in [x_0,b]\right\}\cdot {d\over dz}\P
\left\{S^+_{T^+(s,x_0,\omega)}\in [x_0,z]\right\}dz\\[4mm]
\qquad\ds =~\Theta(x_0) \cdot\P\left\{S^+_{T^+(t+s,x_0,\omega)}\in [x_0,b]\right\}~=~P_{t+s}\bigl(x_0,[x_0,b]\bigr).
\enda $$
The remaining case is when $x_0=y_j\in \S^*$.   We then observe that, for 
every 
$\tau\geq 0$, one has $$P_\tau\bigl(x_0,\{x_0\}\bigr)~=~ \P\bigl\{Y_j(\omega)\geq \tau\bigr\}~=~e^{-\lambda_j\tau},$$ 
\[
P_\tau\bigl(x_0,[x_0,b]\bigr)~
=~\Theta(x_0)\cdot \P\left\{S^+_{T^+(\tau,x_0,\omega)}\in [x_0,b]\right\}+\left(1-\Theta(x_0)\right)\cdot \P\bigl\{Y_j(\omega)\geq \tau\bigr\}.
\]
Using the above identities we obtain
\[\bega{l}
\ds \int P_t\bigl(z,[x_0,b]\bigr) P_s(x_0,dz)~=~P_t\bigl(x_0,[x_0,b]\bigr)\cdot P_s\bigl(x_0,\{x_0\}\bigr)+\int_{]x_0,b]} P_t\bigl(z,[x_0,b]\bigr) P_s(x_0,dz)
\\[4mm]
\qquad\ds =~\bigg[ (1-\Theta(x_0))\cdot \P\{Y_j\geq t\}+\Theta(x_0)\cdot\P\left\{S^+_{T^+(t,x_0,\omega)}\in [x_0,b]\right\}\bigg]\cdot  \P\left\{Y_j\geq s\right\}\\[4mm]
\qquad\qquad  \ds +\Theta(x_0)\cdot\int_{]x_0,b]}\P\left\{S^+_{T^+(t,z,\omega)}\in [x_0,b]\right\}\cdot{d\over dz}\P\left\{S^+_{T^+(s,x_0,\omega)}\in [x_0,z]\right\}dz\\[4mm]
\qquad\ds =~(1-\Theta(x_0))\cdot \P\{Y_j\geq t+s\}
\\[4mm]
\qquad\qquad \ds+\Theta(x_0)\cdot\P\left\{S^+_{T^+(t,x_0,\omega)}\in [x_0,b]\right\}\cdot  \P\left\{S^+_{T^+(s,x_0,\omega)}=x_0\right\}\\[4mm]
\qquad\qquad \ds  +\Theta(x_0)\cdot\int_{]x_0,b]}\P\left\{S^+_{T^+(t,z,\omega)}\in [x_0,b]\right\}\cdot{d\over dz}\P\left\{S^+_{T^+(s,x_0,\omega)}\in [x_0,z]\right\}dz\\[4mm]\ds 
\qquad\ds =~(1-\Theta(x_0))\cdot\P\left\{S^-_{T^-(t+s,x_0,\omega)}=x_0\right\}+\Theta(x_0)\cdot\P\left\{S^+_{T^+(t+s,x_0,\omega)}\in [x_0,b]\right\}\\[4mm]
\qquad \ds =~P_{t+s}\bigl(x_0,[x_0,b]\bigr).
\enda
\]
\endproof
\section{Characterization of Markov semigroups compatible with the ODE}
\label{s:6}
\setcounter{equation}{0}
The goal of the present section is to describe the most general Markov semigroup whose 
sample paths are all solutions to the ODE (\ref{ode}). 
Our main result shows that all of these Markov 
semigroups are of the form considered in 
Proposition~\ref{p:5}.
Namely, they are all obtained starting with a deterministic semigroup and performing two modifications: (i) adding a countable number of random waiting times at points $y_j\in f^{-1}(0)$, each with a Poisson distribution, and 
(ii) assigning the probabilities of moving upwards or downwards, at 
isolated points from where both an increasing and a decreasing solution can initiate.

\begin{theorem}\label{t:61} Let $f$ be a function satisfying {\bf (A1)-(A2)}.
The following statements are equivalent.
\begi
\item[(I)] The random variables $X(t,x_0,\omega)$ yield
a Markov process whose sample
paths are solutions to the ODE (\ref{ode})-(\ref{id}).
\item[(II)] There exist: (i) a positive, atomless Borel measure $\mu$ supported on $f^{-1}(0)$,
(ii) a countable set $\S\subseteq f^{-1}(0)$ of stationary points, (iii) a countable set $\S^*
= \{y_j: j\geq 1\}\subseteq f^{-1}(0)$ and corresponding numbers $\lambda_j>0$
determining the Poisson waiting times,
and (iv) a map $\Theta:\Omega^*\mapsto [0,1]$, such that
the  transition 
kernels $P_t(x_0,A)= \P\{ X(t,x_0,\omega)\in A\}$ coincide with the corresponding ones
constructed 
at (\ref{TP1}), (\ref{TP2}), (\ref{TP3}), (\ref{TP4}).
\endi \end{theorem}
\v
{\bf Proof.} 
The implication (II)$\implies$(I) was proved in Proposition~\ref{p:5}.
Here we need to prove (I)$\implies$(II).
\medskip

{\bf 1.} 
We begin by defining two subsets of $f^{-1}(0)$
\bel{OSD1}
\Omega_X^0~\doteq~\Big\{ x_0\,;~\P\{ X(1, x_0,\omega)= x_0\}=1 \Big\},\eeq
\bel{SD2}
\S^*~\doteq~\Big\{ x_0\,;~0<\P\{ X(1, x_0,\omega)=x_0\} <1 \Big\}.\eeq
Here $\Omega_X^0$ is the set of points where the motion stops forever, while $\S^*$ contains points where the motion stops for a  random time, then starts again.
By the Markov property and the fact that all solutions of (\ref{ode}) are monotone, 
$x_0\in \Omega_X^0$ implies $X(t,x_0, \omega)=x_0$ for all $t\geq 0$ and a.e.~$\omega$. We claim that for each $x_0\in\S^*$, there holds
\bel{SD3} \P\{ X(t, x_0,\omega)=x_0\}~=~e^{-\lambda t}\qquad\forall t\geq 0
\eeq
with $ \lambda\doteq -\log\Big( \P\{ X(1, x_0,\omega)=x_0\}\Big)<\infty$. Indeed, 
for $t\geq 0$, the map 
$$t~\mapsto~ g(t)~=~
\log\Big(\P\{ X(t, x_0,\omega)=x_0\}\Big)~\in ~]-\infty,0]$$ 
is nonincreasing. By the Markov property, it satisfies 
\[
g(0)~=~0,\qquad g(t+s)~=~g(t)+g(s)\qquad\forall t,s>0.
\]
Therefore,  $g(q)=-\lambda\, q$ for some $\lambda>0$ and 
all positive rational numbers $q\in\mathbb{Q}_+$. Since $g$ is monotone, 
we conclude that $g$ is continuous on $[0,+\infty[\,$ and hence 
\[
\P\{ X(t, x_0,\omega)=x_0\}~=~e^{g(t)}~=~e^{-\lambda t}\qquad\forall t\in [0,+\infty[.
\]
%Otherwise, if $ \P\{ X(1, x_0,\omega)=x_0\}=0$ then $\P\{ X(q, x_0,\omega)=x_0\}=0$ for all $q\in \mathbb{Q}^+$. Thus, the decreasing property of the map $[0,+\infty[\ni t\mapsto  \P\{ X(t, x_0,\omega)=x_0\}\in [0,1]$ implies that $\P\{ X(t, x_0,\omega)=x_0\}=0$ for all $t\in [0,+\infty[$.
%
%Next, for $z_k\in \Omega^*\setminus \S$, we define 
%\bel{SD4}\theta_k~=~\Theta(x_k)~\doteq~\P\Big\{ X(t,z_k,\omega)>x_0\quad\hbox{for some}~ t>0\Big\}.\eeq

{\bf 2.} Next, we need to identify
the maximal intervals where the  random
dynamics is  increasing or decreasing.
%\textcolor{blue}{
%\begin{definition}\label{d:rdyn} We say that an  interval $J= [a,b[$ or $J=\,]a,b[$ (open to the right)  is a 
%{\bf domain of increasing random dynamics}
%if, for every $x_0\in J$,  there exists $\Hat x>x_0$ such that 
%\bel{D+}
% \P\Big\{ X(1, x_0, \omega) > \Hat x\Big\}~>~0.\eeq
%Similarly, we say that an  interval $J= ]a,b]$ or $J=\,]a,b[$ (open to the left) is a 
%{\bf domain of decreasing random dynamics}
%if, for every $x_0\in J$,  there exists $\Hat x<x_0$ such that
%\bel{D+}
%\P\Big\{ X(1, x_0, \omega) < \Hat x\Big\}~>~0.\eeq
%\end{definition}
%}
\begin{definition}\label{d:rdynold} We say that an  interval $J= [a,b[$ or $J=\,]a,b[$ (open to the right)  is a 
{\bf domain of increasing random dynamics}
if for every $x_0, \Hat x\in J$ with $x_0<\Hat x$, one has
\bel{D+}
\lim_{t\to \infty} \P\Big\{ X(t, x_0, \omega) > \Hat x\Big\}~=~1.\eeq
Similarly, we say that an  interval $J= ]a,b]$ or $J=\,]a,b[$ (open to the left) is a 
{\bf domain of decreasing random dynamics}
if, for every $x_0, \Hat x\in J$ with $\Hat x < x_0$, one has
\bel{D-}
\lim_{t\to \infty} \P\Big\{ X(t, x_0, \omega) < \Hat x\Big\}~=~1.\eeq
\end{definition}

We observe that, if two intervals of increasing random dynamics have nonempty
intersection, then their union is still a domain of increasing random dynamics.
We can thus define a countable number of maximal open intervals 
$J_k^+=]a_k,b_k[$, $k\in \J^+$ where the random dynamics is increasing, and 
 a countable number of maximal open intervals $J_k^-=]c_k,d_k[$, $k\in \J^-$ where the random dynamics is decreasing. Recalling 
 (\ref{OSD1}), we shall consider the countable set  
 \bel{S-M}
\S~=~\Big(\bigl\{ a_k\,,\,b_k:k\in \J^+\} \cup \bigl\{ c_k\,,\,d_k: k\in \J^-\}
\Big)\cap \Omega_X^0\,. 
\eeq
Next, we consider  the set $\Omega^*$  of isolated points 
from which both a decreasing and an increasing 
trajectory can initiate. For $z_k\in \Omega^*\setminus \S$, we define
the probability of moving upwards as
\bel{SD4}\Theta(z_k)~\doteq~\P\Big\{ X(t,z_k,\omega)>z_k\quad\hbox{for some}~ t>0\Big\}.\eeq
Of course, starting from $z_k$, the probability of moving downwards is thus 
$1-\Theta(z_k)$.
\medskip

{\bf 3.} Consider a maximal domain of increasing random dynamics, say  $J=\,]a,b[$ or $J=[a,b[$. By definition,
 for any $a< x<y <b$ we have that
\bel{pa1}\prob\{ X(\tau, x,\omega) < y\}~\doteq~\delta~<~1\eeq
for some $\tau>0$.
This implies
\bel{pa2}\prob\{ X(k\tau, x,\omega)  < y\}~\leq~\delta^k\qquad\forall k\geq 1.\eeq
Call $T= T^{xy}>0$ the random time needed for a solution starting at $x$ 
to reach $y$. By (\ref{pa2}) it follows
\bel{pa3} \prob\{ T>k\tau\}~\leq~\delta^k\,.\eeq
By (\ref{pa3}), the probability distribution $T=T^{xy}$ has moments of all orders.
In particular, its mean and its variance are finite.

As a consequence for every $x\leq z<z'\leq y$ the random variable 
$T^{zz'}$ has moments of all orders. Its expected value and its variance satisfy
\bel{pa5}
E[T^{zz'}]~\leq~E[T^{xy}], \qquad   \hbox{Var.}(T^{zz'})~\leq~\hbox{Var.}(T^{xy}).
\eeq
Now consider the points $y_j\in [x, y]$ where a random waiting time $Y_j$
occurs.   The Markov property implies that $Y_j$ has a Poisson distribution
as in (\ref{pxj}), with expected value $E\bigl[ Y_j\bigr]=\ds{1\over \lambda_j}$. From (\ref{pa5}), it follows that for every finite subset $S^*_1$ of $S^*\cap [x,y]$
\[
\sum_{y_j\in \S^*_1} {1\over \lambda_j}~=~\sum_{y_j\in \S^*_1}
E\bigl[ Y_j\bigr]~\leq~E[T^{xy}]~<~+\infty.
\]
In particular,
\[
\#\{y_j\in \S^*\cap [x,y]: \lambda_j\leq n\}~\leq~n\cdot E[T^{xy}]\qquad\forall n\geq 1,
\]
and this shows that the set  $\S^*$ of points where a random waiting time occurs is at most countable. Thus, we can write $\S^*\cap J=\{y_1, y_2,\ldots\}$ and 
\bel{fsum}
\sum_{z \leq y_j< z'}{1\over \lambda_j}~=~\sum_{z \leq y_j< z'}E\bigl[ Y_j\bigr]~\leq~ E[T^{zz'}].
\eeq

{\bf 4.}  To help the reader, we provide here an outline the remaining steps of the proof. Let $J=]a,b[$ or $J=[a,b[$ be a maximal domain of increasing random dynamics and  $\{y_1,y_2,\ldots\} = J\cap \S^*$ be the set of points where the trajectories 
 stop for a random time $Y_j$, then start again.    We can then define a new Markov semigroup by removing these waiting times.   More precisely,   for each $n\geq 1$,
 we define a Markov semigroup $S^{(n)}$ whose trajectories $t\mapsto X^{(n)}(t,x_0,\omega)$ satisfy
 $$ X^{(n)} (\tau_n(t,x_0,\omega), x_0, \omega)~=~X(t,x_0,\omega),$$ 
where
$$\tau_n(t,x_0,\omega)~\doteq~t - \meas\Big\{ s\in [0,t]\,;~  X(s, x_0,\omega)\in \{y_1,\ldots,y_n\}\Big\}.$$
In turn, given $S^{(n)}$, we can recover the original semigroup $S$ by inserting back
the waiting times at the points $y_1,\ldots, y_n$.

Now consider the limit semigroup $\Tilde S=\lim_{n\to\infty} S^{(n)}$.
Since we removed all the random waiting times, restricted to the interval 
$J$ the trajectories of $\Tilde S$ are 
strictly increasing with probability one.   This will imply that
$\Tilde S$ is a deterministic semigroup, hence it
admits the representation proved in Theorem~\ref{car_sem}. Namely, there exists
a positive atomless measure $\mu$ on $J\cap f^{-1}(0)$ such that trajectories of $\Tilde S$
are determined by the formula (\ref{St+}).   In turn, the original Markov semigroup $S$
can be recovered from $\Tilde S$ by adding back  the random 
waiting times $Y_j$
at the points $y_j$.   This will provide the desired characterization of $S$.
\medskip

{\bf 5.}
%Since $\S^*$ is countable, we can write $\S^*\cap J=\{y_1, y_2,\ldots\}$ and 
%\bel{fsum}
%\sum_{j=1}^{\infty}{1\over \lambda_j}~=~\sum_{z \leq y_j<z'}
%\eeq
We now work out details. Consider the first point $y_1$.   We claim that, starting
with the Markov semigroup $S$ and removing the random waiting time 
$Y_1(\omega)$ at $y_1$, we obtain another Markov semigroup $S^{(1)}$. 

Indeed, given an initial point $x_0\in J$, let $t\mapsto X(t, x_0,\omega)$
be a sample path of the original process, and define the new path
$X^{(1)}(\cdot, x_0, \omega)$ by setting
\[
X^{(1)}(t,x_0,\omega)~=~
\begin{cases}
X(t,x_0,\omega)&\mathrm{if}\quad X(t,x_0,\omega)\in ]a,y_1[~~\mathrm{or}~~x_0>y_1\,,
\cr\cr
X(t+Y_1(\omega),x_0,\omega)& \mathrm{if}\quad X(t,x_0,\omega)\in [y_1,b[~~\mathrm{and}~~x_0\leq y_1\,.
\end{cases}
\]
We claim that the transition probability kernels 
\[
P^{(1)}_\tau(x_0,A)~\doteq~\P\left\{X^{(1)}(\tau,x_0,\omega)\in A\right\}
\]
satisfy the Chapman-Kolmogorov equation. Indeed, for any Borel set $B\subset ]a,b[\,$, the following holds. 

If $x_0>y_1$, then
\[
P^{(1)}_t(x_0,B)~=~P_t(x_0,B).
\]
Otherwise, if $x_0\leq y_1$, setting $B_1=\, ]a,y_1[\cap B$ and $B_2= [y_1,b[\cap B$,
one obtains
\begin{eqnarray*}
P^{(1)}_t(x_0,B)&=&\P\left\{X(t,x_0,\omega)\in B_1\right\}+\P\bigl\{X(t+Y_1(\omega),x_0,\omega)\in B_2\bigr\}\\
&=&P_t(x_0, B_1)+\ds\int_{0}^{\infty}\P\bigl\{X(t+\tau,x_0,\omega)\in B_2\bigr\}\lambda_1 e^{-\lambda_1 \tau}d\tau\\
&=&P_t(x_0, B_1)+\ds\int_{0}^{\infty} P_{t+\tau}(x_0,B_2)\cdot \lambda_1e^{-\lambda_1 \tau}d\tau\,.
\end{eqnarray*} 
Therefore, for any $s,t>0$ and any Borel set $A\subseteq ]a,b[$,  if $x_0>y_1$ then 
we trivially have
\[
\ds\int P^{(1)}_t(z,A)P^{(1)}_s(x_0,dz)~=~\ds\int P_{t}(z,A)P_s(x_0,dz)~=~P_{t+s}(x_0,A)~=~P^{(1)}_{t+s}(x_0,A).
\]
Otherwise, if $x_0\leq y_1$, then we consider the sets $A_1=A\cap ]a,y_1[\,$, $A_2=A\cap [y_1,b[\,$, and compute
$$
\bega{l}
\ds\int P^{(1)}_t(z,A)P^{(1)}_s(x_0,dz)\\ [4mm]
\qquad\qquad =~\ds\int P_{t}(z,A_1)P_s(x_0,dz)+\ds\int  P_{t}(z,A_2)\cdot \left [\int_{0}^{\infty} P_{s+\tau}(x_0,dz)\cdot \lambda_1e^{-\lambda_1\tau}d\tau\right]\\[4mm]
\qquad\qquad =~ \ds P_{t+s}(x_0,A_1)+\int_{0}^{\infty}\left[\int P_{t}(z,A_2)\cdot  P_{s+\tau}(x_0,dz)\right]\cdot \lambda_1e^{-\lambda_1 \tau}d\tau\\[4mm]
\qquad\qquad =~\ds P_{t+s}(x_0,A_1)+\int_{0}^{\infty}P_{t+s+\tau}(x_0,A_2)\cdot  \lambda_1e^{-\lambda_1 \tau}d\tau~=~P^{(1)}_{t+s}(x_0,A).
\enda
$$
Hence, the process $S^{(1)}$ with sample paths $t\mapsto X^{(1)}(t,x_0,\omega)$ is a Markov semigroup.
\medskip

{\bf 6.} By induction, for each $n\geq 1$, we consider the
process $S^{(n)}$ obtained from $S^{(n-1)}$
by removing the Poisson waiting time at the point $y_n\in J\cap \S^*$.
By the previous argument, every $S^{(n)}$ is a Markov semigroup.  
Conversely, one can recover the original Markov semigroup $S$ starting with 
$S^{(n)}$ and adding the Poisson waiting times $Y_1,\ldots, Y_n$ at the points
$y_1,\ldots,y_n$. If the set $J\cap\S^*$ contains $N$ points, the induction argument 
is concluded  $N$ steps.
Otherwise,  since  the sequence of trajectories $t\mapsto 
X^{(n)}(t,x_0, \omega)$ is monotone increasing and bounded for $t$ in bounded sets, 
for every $\omega$ there exists the limit $X^{(n)}(t,x_0, \omega)\to \Tilde X(t, x_0, \omega)$, uniformly for $t$ in bounded sets.    By part 
(iii) of Theorem~\ref{t:1}, each limit trajectory $t\mapsto \Tilde X(t, x_0, \omega)$
is still a solution of the ODE (\ref{ode}).

In the next step we will prove that this limit process is deterministic.
Namely
\bel{var0} \hbox{Var.}\Big( \Tilde X(t, x_0, \omega)\Big)~=~0,\eeq
for every $x_0\in J$ and $t>0$.  
\v
{\bf 7.} Given $x_0< \Hat x\in J$, let us define 
\[
\begin{cases}
T_n^{x_0 \Hat x}(\omega)&\doteq~\inf\,\Big\{ t>0\,;~~X^{(n)}(t,x_0,\omega)\geq \Hat x
\Big\}\\
\Tilde T^{x_0 \Hat x}(\omega)&\doteq~\inf\,\Big\{ t>0\,;~~\Tilde X(t,x_0,\omega)\geq \Hat x
\Big\}
\end{cases}
\]
the random times needed to reach $\Hat x$,  for a trajectory starting at $x_0$, having removed respectively  the first $n$ waiting times at $y_1,\ldots, y_n$ or all waiting times $y_j\in J\cap S^*$.

Let $\ve_0>0$ be such that $E[\Tilde T^{x_0 \Hat x}]>\ve_0$.
By (\ref{fsum}) we can choose $n$ large enough so that
\bel{var1}
\sum_{j>n,~ x_0\leq y_j< \Hat x}~ {1\over \lambda_j}~<~\ve_0\,.
\eeq

 We now consider
the expected values of these random times.
By construction the map $y \mapsto  E\big[T_n^{x_0 y}\big]$ is nondecreasing, and has upward jumps in the amount $\ds {1\over \lambda_j}$ at each point $y_j$, with $j>n$.  By (\ref{var1}),  we can introduce 
a partition
$$x~=~x_0~<~x_1~<~\cdots~<~x_\nu ~=~\Hat x$$
such that
\bel{dt}\ve_0~\leq~
E[T_n^{x x_i}] - E[T_n^{x x_{i-i}}]~=~E[T_n^{x_{i-1}x_i}]~\leq~2\ve_0\,
\qquad\forall i=1,2,\ldots,\nu.\eeq
Notice that, setting 
$$K~\doteq~ E[T_n^{x_0 \Hat x}],$$
the inequalities (\ref{dt}) require
\bel{nub}
{K\over 2\epsilon_0}~\leq~\nu~\leq~{K\over\epsilon_0}\,.\eeq
From (\ref{dt}) it follows
\[
2\ve_0~\geq~E\left[T_n^{x_{i-1}x_i}\right]~\geq ~{4K\over \nu}\cdot \P\left\{T_n^{x_{i-1}x_i}\geq {4K\over \nu}\right\},
\]
and (\ref{nub}) implies 
\[
\prob\left\{ X^{(n)}\Big({4K\over \nu}, x_{i-1} ,\omega
\Big)~<~x_i\right\}~=~\P\left\{T_n^{x_{i-1}x_i}\geq {4K\over \nu}\right\}~\leq~{1\over 2}\cdot {\nu\ve_0\over K}~\leq~{1\over 2}\,.
\]
In turn, for every $j\geq 1$, it holds 
\begin{eqnarray*}
\P\left\{T_n^{x_{i-1}x_i}\geq {4Kj\over \nu}\right\}&=&\prob\left\{ X^{(n)}\Big({4jK\over \nu}, x_{i-1} ,\omega
\Big)~<~x_i\right\}~\leq~2^{-j}.
\end{eqnarray*}
%Thus, consider a random trajectory starting at $x_{i-1}$. The probability that 
%at time $t = 2K/\nu$ this trajectory has not yet crossed the point $x_i$ is 
%bounded by
%\[
%\prob\left\{ X^{(n)}\Big({4K\over \nu}, x_{i-1} ,\omega
%\Big)~<~x_i\right\}~\leq~{1\over 2}\,.
%\]
%%
%%Because of (\ref{dt}), the probability that 
%%at time $t = 2K/\nu$ this trajectory has not yet crossed the point $x_i$ is 
%%estimated by
%%\bel{pa6}
%%\prob\left\{ X^{(n)}\Big({4K\over \nu}, x_{i-1} ,\omega
%%\Big)~<~x_i\right\}~\leq~{1\over 2}\,.\eeq
%%Indeed, if (\ref{pa6}) fails, then (\ref{nub})  would imply
%%$$E[T_n^{x_{i-1}x_i}]~>~{1\over 2} \cdot {4K\over\nu}~
%%\geq~2\epsilon_0\,,$$
%%in contradiction with (\ref{dt}).   
%In turn, for every $j\geq 1$, it holds 
%\bel{pa7}\prob\left\{ X^{(n)}\Big({4jK\over \nu}, x_{i-1} ,\omega
%\Big)~<~x_i\right\}~\leq~2^{-j}\,.\eeq
Let us now estimate the variance of the random time $\Tilde T^{x_{i-1}x_i}$ needed to move from $x_{i-1}$ to $x_i$. Call 
$$\Phi(s)~\doteq~\prob\{\Tilde T^{x_{i-1}x_i}>s\}.$$
Observing that $\Phi$ is a non-increasing function 
which satisfies
$$\Phi(4jK/\nu)~\leq~\P\left\{T_n^{x_{i-1}x_i}\geq {4Kj\over \nu}\right\}~\leq~2^{-j}\qquad\quad\forall j\geq 1,$$
we estimate the variance\footnote{Indeed, differentiating  the identity
$\ds {1\over 1-x}~=~\sum_{j=0}^\infty x^j$  one obtains
$\ds {1\over (1-x)^2} ~=~\sum_{j=1}^\infty j x^{j-1}$. 
Setting $x=1/2$ we get~ $ \ds\sum_{j=1}^\infty j 2^{1-j}~=~4$.}
\bel{pa8}
\bega{rl}\hbox{Var.} (\Tilde T^{x_{i-1}x_i})& \ds\leq~-\int_0^{+\infty} s^2 d\Phi(s)
~=~\int_0^{+\infty} 2 s\, \Phi(s)\, ds\\[4mm]
&\ds =~\sum_{j\geq 1} \int_{4(j-1)K/\nu}^{4jK/\nu}  2 s\, \Phi(s)\, ds
~\leq~{4K\over \nu} \sum_{j\geq 1} {4jK\over \nu}~2^{1-j}~=~{64\,K^2\over \nu^2}\,.
\enda
\eeq
We  observe that, for $i=1,\ldots,\nu$, 
the random variables $\Tilde T^{x_{i-1} x_i}$ are 
mutually independent.   Hence, by (\ref{nub}),
\bel{pa9}\hbox{Var.} (\Tilde T^{x_0 \Hat x})~=~\sum_{i=1}^\nu \hbox{Var.}( \Tilde T^{x_{i-1} x_i})
~\leq~\sum_{i=1}^\nu {64\,K^2\over \nu^2}~\leq~{64\,K^2\over \nu}~
\leq~128 K \ve_0\,.\eeq
Since  $\ve_0>0$ can be chosen arbitrarily small, we conclude that
the variance of $\Tilde T^{x_0 \Hat x}$ is zero.   Hence,  the limit process is deterministic.
\v
{\bf 8.} We now define
\bel{TST}
\Tilde S_t x_0~\doteq~\lim_{n\to\infty} ~E\Big[ X^{(n)} (t, x_0, \omega)\Big].\eeq
The existence of the limit follows trivially by monotonicity.  
Since the variance of the random variables $X^{(n)}$ goes to zero, again
by (iii) in Theorem~\ref{t:1} it follows that each trajectory $t\mapsto
\Tilde S_t x_0$
is a Carath\'eodory solution to (\ref{ode}).

We claim that $\Tilde S$ is a deterministic semigroup.  Toward this goal, we first observe that the maps
$$t~\mapsto~\Tilde S_t x_0\,,\qquad x_0\mapsto \Tilde S_t x_0$$
are both Lipschitz, and strictly increasing (as long as $\Tilde S_t x_0<b$).
To prove the semigroup identity
\bel{SID}
\Tilde S_s(\Tilde S_t x_0)~=~\Tilde S_{t+s} x_0\,,\eeq
set $x_1 = \Tilde S_t x_0$ and let $\ve>0$ be given.
Choose $\delta>0$ such that 
\bel{dep}\Tilde S_s(x_1-\delta)~>~\Tilde S_s x_1-\ve\,,\qquad\qquad \Tilde S_s(x_1+\delta)~<~
\Tilde S_s x_1 + \ve.\eeq
We now have
$$\lim_{n\to\infty} \prob\Big\{ X^{(n)}(t, x_0, \omega)\in [x_1-\delta, x_1+\delta] \Big\} ~=~1\,.$$
By (\ref{dep}) and a comparison argument, this implies
$$\lim_{n\to\infty} \prob\Big\{ X^{(n)}(t+s, x_0, \omega)\in [\Tilde S_s x_1-\ve, \Tilde S_s x_1+\ve] \Big\} ~=~1\,.$$
Since $\ve>0$ was arbitrary, this proves the semigroup property (\ref{SID}).

As a consequence, by Theorem~\ref{car_sem}, there exists a positive
atomless measure $\mu$ supported on $J\cap f^{-1}(0)$ such that the trajectories
of the semigroup are characterized by the formula (\ref{St+}).

Of course, the above construction of a measure $\mu$ can be repeated on 
every maximal interval $J$ of increasing or decreasing dynamics. 
\medskip

{\bf 9.}   As in Proposition~\ref{p:5}, given  the sets $\S$, $\S^*\subseteq f^{-1}(0)$, the positive, 
atomless measure $\mu$ on $f^{-1}(0)$,  
and the maps $\Lambda:\S^*\mapsto \R_+$, $\Theta:\Omega^*\mapsto [0,1]$,
we construct a unique Markov semigroup $\Hat S$, with transition probabilities given at (\ref{TP1}), (\ref{TP2}), (\ref{TP3}), (\ref{TP4}). 
Recalling Definition \ref{d:increase}, let $I_i^+$, $i\in \I^+$, and $I_i^-$, $i\in \I^-$, be respectively the domains of increase and the domains of decrease corresponding to $\mu$ and $\S$.  For any $k\in \J^+$ and $x_0\in J^+_k$ with $J_k^+=]a_k,b_k[$ or $J_k^+=[a_k,b_k[$, the map $[0,T_{b_k}[\ni t\mapsto \Tilde{S}_tx_0$ is strictly increasing and $\lim_{t\to T_{b_k}}\Tilde{S}_tx_0=b_k$ for some $T_{b_k}>0$. Thus, $J^+_k$ is an strictly increasing domain associated to $\Tilde{S}$ and this implies that 
$$
J^+_k~\subseteq~ I^+_i\qquad\mathrm{for~some}~i\in \I^+.
$$
Similarly, for every $h\in\mathcal{J}^-$ it holds that $J^-_h~\subseteq~ I^-_j$ for some $i\in\mathcal{I}^-$. Thus, 
$$
\bigcup_{k\in\J^+}J^+_k~\subseteq~\bigcup_{i\in \I^+} I_i^+\qquad\mathrm{and}\qquad \bigcup_{h\in\J^-}J^-_h~\subseteq~\bigcup_{i\in \I^-} I_i^-.
$$

In the remainder of the proof, we will show that the given Markov semigroup $S$ coincides with the semigroup $\Hat S$ constructed in Proposition~\ref{p:5}.
In other words, for every $x_0\in\R$, $t>0$ and any Borel set $A\subset\R$, the transition probabilities coincide:
\bel{PSH}
P_t(x_0,A)~=~\Hat P_t(x_0,A)~\doteq~\P\bigl\{ \Hat S_t x_0\in A\bigr\}.\eeq
Four cases must be considered:
\begin{itemize}
\item [(i)]  $x_0\notin\left(\bigcup_{i\in\I^+} I_i^+\right)
\cup\left(\bigcup_{i\in \I^-} I_i^-\right)$;
\item [(ii)] $x_0\in \left(\bigcup_{i\in \I^+} I_i^+\right)\setminus
\left(\bigcup_{i\in\I^-} I_i^-\right)$;
\item [(iii)] $x_0\in \left(\bigcup_{i\in \I^-} I_i^-\right)\setminus
\left(\bigcup_{i\in\I^+} I_i^+\right)$;
\item[ (iv)] $x_0\in \left(\bigcup_{i\in I^-} I_i^-\right)\cap
\left(\bigcup_{i\in \I^+} I_i^+\right)\subseteq\Omega^*$.
\end{itemize}

Case (i) is trivial. We  thus focus on Case (ii). Assume that $x_0\in I_i^+\setminus
\left(\bigcup_j I_j^-\right)$ for some $i\in\I^+$.  In this case, we only need to prove (\ref{PSH}) for $A= [x_0, \Tilde{S}_s(x_0)]$ with $s>0$. From the previous construction in Proposition~\ref{p:5}, it holds that 
\bel{TPH2}
\Hat P_t\bigl( x_0,\, [x_0, \Tilde{S}_s(x_0)]\bigr)~=~\P\left\{ s+\sum_{y_j\in S^*\cap [ x_0, \Tilde{S}_s(x_0)]}  Y_j (\omega)~\geq ~t\right\}.
\eeq
Two cases are considered: 
%\bel{TPH2}\bega{l}
%\ds\Hat P_t\bigl( x_0,\, \,]-\infty, x_0[\bigr) ~=~0,\\[4mm]
%\ds \Hat P_t\bigl( x_0,\, [x_0, \Tilde{S}_s(x_0)]\bigr)~=~\P\left\{ s+\sum_{y_j\in S^*\cap [ x_0, \Tilde{S}_s(x_0)]}  Y_j (\omega)~\geq ~t\right\}.\enda
%\eeq
\begin{itemize}
\item  If $x_0\in J_k^+$ for some $k\in \J^+$ then one  has
\[\bega{l}
\ds P_t\left(x_0, [x_0,\Tilde{S}_s(x_0)]\right)~=~\P\left\{X(t,x_0,\omega)\in [x_0,\Tilde{S}_s(x_0)]\right\}\\[4mm]
\qquad =~\ds\P\left(\Tilde{X}(\tau(t,x_0,\omega),x_0,\omega)\in  [x_0,\Tilde{X}(s,x_0,\omega)]\right)\\[4mm]
\qquad\ds=~\P\left(\tau(t,x_0,\omega)\leq s\right)~=~\P\left\{ s+\sum_{y_j\in S^*\cap [ x_0, \Tilde{S}_s(x_0)]}  Y_j (\omega)~\geq ~t\right\},\enda
\]
where
$$\tau(t,x_0,\omega)~\doteq~t - \meas\Big\{ r\in [0,t]\,;~  X(r, x_0,\omega)\in \S^*\Big\}.$$
This yields (\ref{PSH}) for  $A= [x_0, \Tilde{S}_s(x_0)]$.
\item Otherwise, if $x_0\notin \bigcup_{k\in \J^+}J^+_k$ then 
\bel{atx}
\lim_{t\to\infty} \P\{X(t,x_0,\omega)\leq x\}~=~1\qquad\forall x>x_0. 
\eeq
and this implies that 
\bel{11}
\P\bigl\{ S_t x_0=x_0\bigr\}~=~1\qquad\mathrm{and}\qquad x_0~\in~\Omega^0_X.
\eeq
Indeed, assume by a contradiction that there exists $\Hat x>x_0$ such that $\P\{X(\Hat t,x_0,\omega)\leq \Hat x\}=\delta<1$ for some $\Hat t>0$. By the Markov property of $X$, we have 
\[
\P\{X(k\Hat t,x_0,\omega)\leq \Hat x\}~\leq~ \delta^k\qquad\forall k\in\mathbb{Z}^+.
\]
In particular, the  monotone increasing property of the maps $t\mapsto X(t,x_0,\omega)$ implies that  $\lim_{t\to\infty}\P\{X(\Hat t,x_0,\omega)\leq \Hat x\}=0$. Hence, $x_0$ is in a domain of increasing random dynamics and this  yields a contradiction. 
\medskip

On the other hand, since $x_0\notin \left(\bigcup_{k\in \J^+}J^+_k\right)\bigcup\left(\bigcup_j I_j^-\right)$ and $\bigcup_{h\in\J^-}J^-_h\subseteq\bigcup_{i\in \I^-} I_i^-$, one has that $x_0\in \Big(\bigl\{ a_k\,,\,b_k:k\in \J^+\} \cup \bigl\{ c_k\,,\,d_k: k\in \J^-\}\Big)$. Thus, (\ref{11}) and (\ref{TP1}) implies that $x_0\in\S$ and 
\[
\Hat P_t\bigl( x_0, \{x_0\}\bigr)~=~1~=~P_t\bigl( x_0, \{x_0\}\bigr).
\]
\end{itemize}
Similarly, one can prove (\ref{PSH}) for $x_0$ in an interval of decreasing dynamics. 

Finally, consider $x_0\in\left(I_i^+\cap I_{i'}^-\right)\backslash \S$, for some $i\in \I^+$ and $i'\in\I^-$. It is sufficient to verify (\ref{PSH})  for $A=[x_0,b]$, with $b\in I_i^+$. The case $A=[c,x_0]$, $c\in I_{i'}^-$ is entirely similar. If $x_0\notin \S^*$ then  
$$
\P\left\{X(1,x_0,\omega)=x_0\right\}~=~0\qquad\mathrm{and}\qquad \P\{X(1,x_0,\omega)>x_0\}~=~\Theta(x_0)
$$
and this implies
$$\bega{l}
\ds \P\left\{X(t,x_0,\omega)\in [x_0,b]\right\}~=~\P\left\{X(t,x_0,\omega)\in [x_0,b]~~\hbox{and}~~X(1,x_0,\omega)>x_0\right\}\\[3mm]
\qquad =~\ds\P\{X(1,x_0,\omega)>x_0\}\cdot \P\left\{X(t,x_0,\omega)\in [x_0,b]\big| X(1,x_0,\omega)>x_0\right\}\\[3mm]
\ds \qquad =~\Theta(x_0)\cdot \P\bigl\{\Tilde S^+_{T^+(t, x_0, \omega)}\in [x_0,b]\bigr\}~=~\Hat P_t(x_0,[x_0,b]).
\enda $$
Otherwise, if $x_0=y_j\in \S^*$ then 
\[
\bega{l}
\P\bigl\{X(t,x_0,\omega)\in [x_0,b]\bigr\}\\[3mm]
\qquad \ds=~\P\Big\{X(t,x_0,\omega)\in [x_0,b]~~\hbox{and}~~X(s,x_0,\omega)>x_0~\mathrm{for~some}~s>0\Big\}\\[3mm]
\qquad \qquad +\P\Big\{X(t,x_0,\omega)\in [x_0,b]~~\hbox{and}~~X(s,x_0,\omega)<x_0~\mathrm{for~some}~s>0\Big\}\\[3mm]

\qquad =~\Theta(x_0)\cdot  \P\bigl\{\Tilde S^+_{T^+(t, x_0, \omega)}\in [x_0,b]\bigr\}\bigg\}\\[3mm]
\qquad\qquad +(1-\Theta(x_0))\cdot \P\Big\{X(t,x_0,\omega)=x_0\,\Big| ~X(s,x_0,\omega)<x_0~\mathrm{for~some}~s>0\Big\}\\[3mm]
\qquad=~\Theta(x_0)\cdot \P\left\{S^+_{T^+(t,x_0,\omega)}\in [x_0,b]\right\}+(1-\Theta(x_0))\cdot  \P\bigl\{Y_j(\omega)\geq t\bigr\}\\[3mm]
\qquad =~\Hat P_t\bigl(x_0,[x_0,b]\bigr).
\enda
\]
This establishes the identity (\ref{PSH}) for all initial points $x_0\in\R$, completing the proof of the theorem.
\endproof

{\bf Acknowledgments.} This research by K. T. Nguyen was partially supported by a grant from
the Simons Foundation/SFARI (521811, NTK).

\vs

\end{document}